\documentclass{article}
\usepackage{graphicx} 
\usepackage{amsthm}
\usepackage{amsmath}
\usepackage{amsfonts}
\newtheorem{proposition}{Proposition}
\newtheorem{lemma}{Lemma}
\newtheorem{theorem}{Theorem}

\newtheorem{corollary}{Corollary}
\usepackage[colorlinks]{hyperref}
\usepackage{cleveref}
\usepackage{xcolor}
\usepackage{cancel}
\usepackage[normalem]{ulem}
\usepackage{authblk}
\newtheoremstyle{boldremark}
  {3pt} 
  {3pt} 
  {\upshape} 
  {} 
  {\bfseries} 
  {.}
  { }
  {} 

\theoremstyle{boldremark}
\newtheorem{remark}{Remark}

\newcommand{\R}{\mathbb{R}}
\newcommand{\N}{\mathbb{N}}
\newcommand{\Sph}{\mathbb{S}}
\newcommand{\Fc}{\mathcal{F}}
\newcommand{\Pp}{\mathbb{P}}
\newcommand{\E}{\mathbb{E}\,}
\newcommand{\Poi}{\mathrm{Poisson}}
\newcommand{\eps}{\varepsilon}
\newcommand{\ONE}{\mathbf{1}}
\newcommand{\sgn}{\mathop{\mathrm{sgn}}}
\newcommand{\as}{\mathrm{a.s.}}
\newcommand{\aseq}{\stackrel{\as}{=}}

\newcommand{\Dwv}{\delta}
\newcommand{\dd}{\mathrm{d}}

\title{Differentiability of Shape Function for Piecewise Smooth Lagrangians in Poissonian Random Environments.}

\author[1]{Yuri Bakhtin\thanks{bakhtin@cims.nyu.edu}} 

\author[1]{Benjamin Edidin\thanks{be2282@nyu.edu}}
\affil[1]{Courant Institute of Mathematical Sciences, New York University}
\date{}

\makeatletter
\newcommand{\namedlabel}[2]{\def\@currentlabel{#2}\label{#1}}
\makeatother

\begin{document}

\maketitle
\begin{abstract}
    We consider a family of models of directed first passage percolation type in continuous multi-dimensional Euclidean space and discrete time with admissible paths constrained to visiting a Poissonian point at each step. We show that even if the Lagrangian defining the cost or energy of a step is allowed to have gradient discontinuities, the resulting shape function, or effective Lagrangian, belongs to the class $C^1$ and provide a formula for its gradient. We also study the limit shapes for the associated growth model. We show that the part of the boundary of the limit shape that can be called the effective front is also $C^1$, and we describe the complementary part, the footprint.  We include a study of  several concrete examples and several simulations.
\end{abstract}

\noindent\textbf{Keywords:} HJB equations, random environment, shape function, limit shape, homogenization, effective Lagrangian, shear invariance, directed polymers, First / Last Passage Percolation

\smallskip

\noindent\textbf{MSC 2020:} 60K37, 82B44, 35B27, 49L12

\section{Introduction}
\label{sec:Intro}
For many variational problems in random environments, the 
growth of the random minimal cost of a path connecting two points grows linearly in
the distance between these points, although the deterministic rate of linear growth depends on the direction. This rate, as a function of direction, is usually called the shape function. In the context of homogenization for stochastic 
Hamilton--Jacobi--Bellman (HJB) equations,
it is called the effective Lagrangian.

Shape functions are always convex but their further regularity properties are generally unknown. For a handful of models possessing special structural properties (symmetries, exact solvability), 
shape functions have been computed explicitly or up to some constants.
For models in the KPZ universality class, it is believed that the associated shape functions are strictly convex and differentiable. These properties may be violated outside the KPZ universality class, see   ~\cite{Bakhtin-Dow:MR4945174},~\cite{bakhtin2023differentiability},~\cite{BakhtinDow-AIHP1529},~\cite{BK18} for discussion and references. 

A general strict convexity result is still one of the central open problems in disordered systems but several differentiability results have been established recently in~\cite{Bakhtin-Dow:MR4945174},~\cite{bakhtin2023differentiability},~\cite{BakhtinDow-AIHP1529}. 
In a recent preprint~\cite{OAI:Strict-convexity-and-differentiability-of-the-planar-exponential-first-passage-limit-shape-September-24-2026},
for First Passage Percolation (FPP) 
with i.i.d.\ edge weights on~$\mathbb{Z}^2$,  differentiability of the limit shape 
for Gamma weights and strict convexity of the limit shape for 
exponential weights are claimed.

The papers~\cite{Bakhtin-Dow:MR4945174},~\cite{bakhtin2023differentiability},~\cite{BakhtinDow-AIHP1529} describe various classes of situations where there is no symmetry or any exact solvability, the shape function is unknown,  but its  differentiability property is rigorously established. These include homogenization for certain stochastic HJB equations, directed and unconstrained FPP,  and directed polymers in continuous space, growth in random Riemannian geometries, etc.

In these papers, one of the important ingredients was the $C^2$ assumption on the Lagrangian of the variational problem. Seemingly central to the argument, this assumption excludes multiple interesting models.  In the present paper, we prove differentiability of the shape function for a model where this~$C^2$~property can be substantially relaxed. In particular, our results imply that 
{\it piecewise}~$C^{1+\alpha}$ Lagrangians with multiple {\it gradient discontinuities} still result in differentiable shape functions. 

For any $d\in\N$, we consider the Poisson process of unit intensity on $\N\times\R^d$. We also equip Poissonian points $(k,x)$ with i.i.d. weights $F(k,x)=F_\omega(k,x)\ge0$ with common distribution with finite expectation.

An admissible path $\gamma:\{0,1,2\ldots,n\}\to\R^d$ is required to start at $0$ and pass only through Poissonian points at intermediate 
times, i.e., for all $k=1,\ldots, n-1$, we require $(k,\gamma_k)$ to be one of the Poissonian points. The action/energy/cost of the path is then given by
\begin{equation}
\label{eq:action}
A^n_\omega(\gamma)=\sum_{k=0}^{n-1}L(\gamma_{k+1}-\gamma_k)+\sum_{k=1}^{n-1} F_\omega(k,\gamma_k).
\end{equation}
where a Lagrangian $L:\R^d\to\R$ is a function satisfying several assumptions that we will state below.  We stress the dependence of the action on the random realization of the Poisson points and weights denoted by $\omega$. 
Such energies first appeared in the literature in~\cite{Peled:MR4586218}, in the context of random surfaces under Gibbs distributions. 

We are interested then in the shape function defined by
\begin{equation}
    \label{eq:shape-f}
    \Lambda(v)=\lim_{n\to\infty}\frac{A^n_\omega(nv)}{n}.
\end{equation}
Here
\begin{equation}
    \label{eq:action-p2p}
    A^n_\omega(x)=\inf \big\{A^n_\omega(\gamma):\gamma\in\Gamma^n_\omega,\ \gamma_n=x\big\},
\end{equation}
where $\Gamma^n_\omega$ is the random collection of admissible paths (starting at $0$ and passing through Poissonian points).

Variational problems of this kind are central to the study of solutions of HJB equations and $\Lambda$ is known as the effective Lagrangian in this context. Thus the results of this paper can be interpreted in terms of HJB homogenization.

If $L(v)=|v|^2$, then this model is shear-invariant and can be studied with the methods of \cite{BCK:MR3110798}, \cite{kickb:bakhtin2016}, \cite{Bakhtin-Li:MR3911894}. In particular, in this case  the shape function is given by $\Lambda(v)=C+v^2$ for some $C$ (we recall this derivation in Section~\ref{sec:Examples}). The differentiability and strict convexity of $\Lambda$ obviously hold.

In~\cite{BakhtinDow-AIHP1529}, a model similar to the described above was considered except the paths were not restricted to Poissonian points and picked up contributions from a random continuous environment.  If one requires $L$ of our model~\eqref{eq:action} to satisfy the (multi-dimensional version of) conditions from~\cite{BakhtinDow-AIHP1529} including the requirement that $L\in C^2$, then one can apply the approach from that paper with minimal adjustments and also establish the differentiability of $\Lambda$.

\medskip

The goal of this paper is to introduce much broader assumptions on $L$ that still guarantee the differentiability of $\Lambda$ and provide a formula for its gradient, see our Theorem~\ref{thm: main}.

We are mostly motivated by including functions $L$ that are piecewise smooth and have discontinuities of
the gradient along a collection of hypersurfaces.
A natural class of functions $L$ not differentiable at $0$ is given by norms and some of their powers. Norms like the $L^1$- and $L^\infty$-norm have whole surfaces of gradient discontinuity providing examples that can be treated by our new results but not by the existing results of \cite{Bakhtin-Dow:MR4945174},~\cite{bakhtin2023differentiability},~\cite{BakhtinDow-AIHP1529}.

One notoriously simple example of this kind is $L(v)=|v|$ in $d=1$.  The formula for the gradient of the associated shape functions takes an extremely simple form. 
We hope to use this formula to prove the shape function's strong convexity near zero for this model,  
see the discussion in Section~\ref{sec:Examples}. It is worth mentioning that in 
~\cite{Hassler:MR5100646} a special lattice FPP model is interpreted in terms of an action minimizing problem with $L(v)=|v|$, and a modification of the approach from   
~\cite{Bakhtin-Dow:MR4945174},~\cite{bakhtin2023differentiability},~\cite{BakhtinDow-AIHP1529} is used to analyze strict convexity of the limit shape at a special point corresponding to infinite speed.

Some other examples are given by functions combining
various norms, e.g.,
\begin{equation}
\label{eq:abs v1 plus squared v2}
L(v)=|v_1|+|v_2|^2,\quad v\in\R^2,    
\end{equation}

or various functions with explicitly piecewise definitions.

\medskip 

A useful point of view at the linear growth of cost typical for studies of models of the First Passage Percolation (FPP) type is to interpret it as 
a growth model and 
look at 
random regions~$E(r)$ consisting of points that can be reached from the origin at a cost not exceeding~$r>0$.
 These regions grow approximately linearly with $r$, and
$\frac{1}{r}E(r)$ converges to a certain limit shape as $r\to\infty$. In our setting, the new result on differentiability of the shape function $\Lambda$ translates into a new result on differentiability of the part of the boundary of the limit shape that can be viewed as an effective front of growth or signal propagation, see Theorems~\ref{thm: limit shape},~\ref{thm: limit shape C1}. The complementary part of the boundary of the limit shape, its footprint, is described in Theorem~\ref{thm:footprint}.

\medskip

The plan of the paper is as follows. In Section~\ref{sec:Main Results}, we describe the setting precisely and give rigorous statements of results. 
In Section~\ref{sec:Examples}, we discuss several examples and provide several plots.
In Section~\ref{sec:main_proofs}, we give a proof of the main result on the shape function differentiability. In Section~\ref{sec:limit-shape-proofs}, we give proofs of results on limit shapes and effective fronts.

{\bf Acknowledgments.} We are grateful to NSF for
partial support through award DMS-2246704.  At the final stage of preparing this paper, Anthropic's Claude was used to check the arguments and the exposition. This helped us avoid several sloppy statements and improve the presentation.

\section{The Setting And Main Results}
\label{sec:Main Results}

In this paper, we work with a model in $1+d$ dimensions with arbitrary spatial dimension $d\in\N$. 
Paths in  
\[
\Gamma^n=\big\{\gamma:\{0,1,\ldots,n\}\to\R^d\,|\, \gamma_0=0 \big\},\quad n\in\N,
\]
may be viewed as directed polymers, i.e., directed chains of monomers interacting with each other and the random environment. They can also be viewed as $d$-dimensional trajectories of Lagrangian particles 
accelerated by a dynamic random environment.

We are interested in variational problems for random energy functionals on paths from $\Gamma^n$.  
For $n\in\N$, every path $\gamma\in\Gamma^n$ is assigned a random energy/ action / cost.
\begin{equation}
    A^n_\omega(\gamma) = \sum_{k=0}^{n-1} L(\Delta_k \gamma)+\sum_{k=1}^{n-1} F_\omega(k,\gamma_k),\quad \omega\in\Omega,
\end{equation}
where 
$
\Delta_k\gamma = \gamma_{k+1} - \gamma_k.
$
The function $L:\R^d\to\R$ is the energy of the nearest neighbor interaction (or the kinetic Lagrangian in the dynamical particle interpretation) and $F_\omega:\N\times\R^d\to\R\cup\{+\infty\}$ is the random potential describing the interaction of the monomers with the environment. Our requirements on these functions and a description of the probability space
$(\Omega,\Fc,\Pp)$ will be given below.

We will assume that $F_\omega=+\infty$ everywhere except the collection of points forming a Poisson point process (PPP) driven by the Lebesgue measure on the spacetime $\N\times\R^d$. Equivalently, layers 
$\{k\}\times\R^d$, $k\in\N$, support mutually independent PPPs driven by the Lebesgue measure (volume) $m(\cdot)$ on $\R^d$. In addition, all Poisson points $(k,x)\in\N\times\R^d$ are equipped with i.i.d.\ strength labels $\varphi(k,x)\ge0$, and we will informally set
\begin{equation}
    \label{eq:potenital}
F_\omega(k,x)=\begin{cases}
    \varphi(k,x),& (k,x)\ \textrm{is a Poisson point,}\\
    +\infty, &\ \textrm{otherwise.}
\end{cases}
\end{equation}

A more rigorous way to do this is to introduce a Borel probability measure $\nu$ on~$\R_+$ playing the role 
of the common distribution of r.v.'s $\varphi(k,x)$ and 
introduce a PPP on the space $\N\times\R^d\times \R_+$ driven by
the measure $\mu=\sum_{k\in\N}\delta_k\otimes m\otimes \nu$. More precisely, 
to define the probability space $(\Omega',\mathcal F',\Pp')$, we first  
define $\Omega'$ as the set  of
locally finite configurations of points $\omega\subset \N\times\R^d\times \R_+$. These point configurations can be identified with $\N\cup\{0\}$-valued measures defined by $\omega(C):=\#(\omega \cap C)$ for bounded Borel sets $C\subset\N\times\R^d\times \R_+$. We define $\Fc'$ to be the $\sigma$-algebra generated by the maps $\omega \mapsto \omega(C)$ for all such sets. The probability measure $\Pp$ on $\Fc$ is uniquely defined by the requirement that under~$\Pp$, for every bounded Borel $C$,  
$\omega(C)\stackrel{d}{=} \Poi(\mu(C))$,
and $\omega(C_1),\dots,\omega(C_n)$ are independent whenever $C_1,\dots,C_n$ are disjoint.

We define the projection $\pi:\N\times\R^d\times \R_+\to \N\times\R^d$ by $\pi(k,x,\varphi)=(k,x)$.
For every Poissonian point $(k,x,\varphi)\in\omega$, the coordinate $\varphi\in\R_+$ can be viewed as a label assigned to
its spacetime footprint $(k,x)$. The collection $\pi\omega$ of these footprints forms a PPP on $\N\times\R^d$ with driving measure $\sum_{k\in\N}\delta_k\otimes m$. 

Let $\Omega$ be the set of configurations $\omega$ such that the spacetime footprints of distinct points of $\omega$ are distinct. Since $\Pp'(\Omega)=1$, from now on we can work on the probability space $(\Omega,\mathcal F,\Pp)=(\Omega,\mathcal F|_\Omega,\Pp'|_\Omega)$. 
On $\Omega$, we can rewrite~\eqref{eq:potenital} more formally.
For every $\omega\in\Omega$,
we can define the map $F_\omega:\N\times\R^d \to \R_+\cup\{\infty\}$ such that if
 $(k,x)\in \pi\omega$, then there is a unique label $F_\omega(k,x)\in\R_+$ such that 
$(k,x,F_\omega(k,x))\in\omega$, and $F_\omega(k,x)=+\infty$ otherwise. 

The only assumption we need to make on this part of the setup is that 
 the i.i.d.\ labels have a finite first moment.

\begin{equation}
\label{eq:labels are L1}
\int_{\R_+}\, \varphi\, \nu(d\varphi)<\infty. 
\end{equation}
We will refer to this part of the setup including condition~\eqref{eq:labels are L1} as
{\bf Assumption~(F)\namedlabel{assum:F}{\rm (F)}.}

We also need to introduce several sets of assumptions on the kinetic Lagrangian $L$. 
The idea is to replace the $C^2$ assumption from~\cite{Bakhtin-Dow:MR4945174}, \cite{bakhtin2023differentiability}, and \cite{BakhtinDow-AIHP1529} by a much less restrictive setup allowing for Lagrangians whose gradients have only H\"older regularity and large sets of discontinuity.

First, we introduce some notation. The standard inner product in $\R^d$ is denoted by $\langle \cdot,\cdot\rangle$, and the Euclidean norm is denoted by $|x|=\langle x,x\rangle^{1/2}$. 
For $x\in\R^d$ and $S\subset\R^d$, we denote $d(x,S)=\inf \{|x-y|:\ y\in S\}$.
For $S\subset \R^d$ and $\eps>0$, we denote $S_{\eps}:= \{x \in \R^d: d(x, S) < \eps \}$.
Euclidean balls are denoted by $B(x,R)=\{y\in\R^d: |y-x|<R\}$ for $x\in\R^d$ and $R>0$.
The volume of the unit ball in $\R^d$ is denoted by  $V_d=\frac{\pi^{d/2}}{\Gamma\left(\frac{d}{2}+1\right)}$.

Assumption~\ref{assum:L1} below consists of several requirements on $L$. We need them
to ensure the existence of the shape function, see Theorem~\ref{thm: shape}.

\smallskip
\textbf{Assumption (L1)\namedlabel{assum:L1}{\rm (L1)}.}
\begin{itemize}
\item
There are numbers $a>0$ and $b\in\R$ such that 
\begin{equation}\label{eq: growth of L}
L(v)\ge a|v|+b,\quad v\in\R^d.
\end{equation}
\item
There are positive constants $C_1$ and $C_2$
such that 
\begin{equation}
\label{eq:comparing to ball volume}
C_22^{d-1}<V_d
\end{equation}
and
\begin{equation}
    \label{eq:assumption on growth of L}
    L(v)\le C_1e^{C_2|v|^d},\quad v\in\R^d.
\end{equation}
\item There is $c>0$ satisfying
\begin{equation}
\label{eq:parameter c}
c^dV_d>1,    
\end{equation}
and such that
\begin{equation}
\label{eq:cond on large increments of L}
\rho_c(s)=o(s),\quad s\to\infty,
\end{equation}
where
\begin{equation}
\label{eq:rho}
\rho_c(s)=\sup\Big\{L(w)-L(v):\ L(v)\le s,\ |w-v|\le c(\ln s)^{1/d}\Big\}.
\end{equation}

\end{itemize}
These conditions hold for a broad class of piecewise smooth $L$ growing polynomially and such that their (weak) derivatives grow polynomially at a slower rate. We will give a more detailed discussion once we introduce further conditions on~$L$. Here we only note that a weaker condition $\lim_{|v|\to\infty}L(v)=+\infty$  was used instead of~\eqref{eq: growth of L} in~\cite{BakhtinDow-AIHP1529}. These conditions are equivalent for convex $L$ but we stress that we do not require $L$ to be convex.

For any $n\in\N$ and $x\in \mathbb{R}^d$, we define
\[
\Gamma^n(x) := \{ \gamma \in \Gamma^n: \gamma_n = x\}.
\]

We are interested in the minimization problem 
\begin{equation}\label{eqn:Anx}
    A^n_\omega(x) = \text{inf}\{A^n_\omega(\gamma) : \gamma \in \Gamma^n(x) \},\quad (n,x)\in \N\times\R^d,
\end{equation}
and in regularity of the linear growth (as $n\to\infty$) of $A^n(nv)$ as a function of direction~$v\in\R^d$.

\begin{theorem}\label{thm: shape} Under Assumptions~\ref{assum:F} and~\ref{assum:L1},
there is a deterministic convex function $\Lambda : \mathbb{R}^d \rightarrow \mathbb{R}$ such that 
\begin{equation}
    \Lambda(v) \aseq \lim_{n\to\infty}\frac{1}{n}A^n_\omega(nv), \quad v\in\R^d.
\end{equation}
\end{theorem}

Although the main idea behind the proof of this theorem is the standard subadditivity argument similar to that in \cite{BakhtinDow-AIHP1529}, the optimal action $A^n$ is actually not exactly subadditive. Thus there is a subtlety,  and we provide a proof in Section~\ref{sec:conv-to-shape}, where we use Assumption~~\ref{assum:L1} to control deviations from a truly subadditive quantity.

\medskip
The limit $\Lambda$ characterizing linear growth as a function of direction is called the shape function. 

To state our main result on differentiability of shape function, we need extra assumptions on $L$:

\textbf{Assumption (L2)\namedlabel{assum:L2}{\rm (L2)}.} There are functions $H, K: \R^d \to \R_+,$ $M: \R^d \to \R^{d}$, a set $S \subset \R^d,$  vectors $u_1, ..., u_{d+1}\in\R^d$ 
and numbers $\alpha,\eps_0, r_0>0$ satisfying the following conditions:
\begin{itemize}
\item
 The interior of the convex hull of $u_1, ..., u_{d+1}$ contains $0$.
\item 
For all $v\in\R^d$ , all $k=1,2,\dots  d+1,$ and all $\eps\in(0,\eps_0) $, 
\begin{equation}\label{eq: linear}
    L(v + \eps u_k) \leq L(v) +  \langle M(v), \eps u_k\rangle + \eps^{1+ \alpha}H(v)\ONE_{v \notin S_\eps} + \eps K(v) \ONE_{v \in S_\eps}.
\end{equation}
\item  $M$ is locally bounded. 
\item There are some $c_0, c_1> 0$ and  $\theta\in[0,1)$ such that for all $v\in\R^d,$ 
\begin{equation}
\label{eq: H}
     H(v)\ONE_{v\notin S_{r_0}}\le c_0+c_1|L(v)|
\end{equation}
and
\begin{equation}\label{eq: K}
     K(v)\le c_0 + c_1|L(v)| ^\theta .
\end{equation}
\item There exists some $c_2 > 0$ such that for all $\eps\in(0,\eps_0)$
\begin{equation}\label{eq: H local}
    \eps ^ \alpha H(v) \ONE_{v\in S_{r_0}\setminus S_\eps} \leq c_2.
\end{equation}

\item For all $R>0$, there exists $\eta_R, k_R > 0$ such that for sufficiently small $\eps>0$,
\begin{equation}\label{eq: sausage}  
m\big( B(0, R) \cap S_\eps \big) \leq k_R \eps^{\eta_R}.
\end{equation}

\end{itemize}

It is easy to see that these conditions imply that $K$ is locally bounded and $L$ is continuous. If one imposes conditions on $L$ from~\cite{BakhtinDow-AIHP1529}, i.e.,  $L\in C^2$ and $\limsup_{|v|\to\infty} (|\nabla^2 L(v)|/L(v))<\infty$, then the required upper bound~\eqref{eq: linear} (with $\eps_0=1$, $\alpha=1$, $K\equiv0$, $M(v)=\nabla L(v)$, and $H(v)=\max_{u:|u-v|\le 1}|\nabla^2 L(u)|$) follows from the second order Taylor expansion. 
The most prominent class of functions satisfying our requirements consists of functions that have well-defined H\"older gradients when restricted to any of several domains that form a partition of $\R^d$, but on the boundaries of these domains the gradient has a jump discontinuity.
The set $S$ in the assumptions above plays the role of the set of gradient discontinuities. 
If not for this set, the proof would only be a slight modification of the proof in~\cite{BakhtinDow-AIHP1529}. Near this set, the Taylor expansion fails, and the term $\eps K(v) \ONE_{v \in S_\eps}$ in~\eqref{eq: linear} is an obstacle for a direct implementation of the approach 
of~\cite{BakhtinDow-AIHP1529}. 
Condition \eqref{eq: sausage} addresses this issue. It means that the set of these discontinuities is not too large.

If Assumption~\ref{assum:L2} holds for $L$ admitting a polynomial upper bound and  a lower bound~\eqref{eq: growth of L} and the function $M$ from \ref{assum:L2}  satisfies  $|M(v)|\le c_0+c_1|L(v)|^\theta$ for some $c_0,c_1>0,\theta\in[0,1)$ and all $v\in\R^d$, then Assumption~\ref{assum:L1} also holds. This defines a large class of functions. Several interesting examples of functions~$L$ satisfying these conditions (but not the conditions of \cite{BakhtinDow-AIHP1529}) are given in Section~\ref{sec:Examples}.

Now we can state our main result.
\begin{theorem}\label{thm: main}
    Suppose Assumptions~\ref{assum:F},~\ref{assum:L1}, and~\ref{assum:L2} hold. Then the shape function~$\Lambda$ defined in Theorem~\ref{thm: shape} belongs to the class $C^1(\R^d)$. For every $v\in \R^d$, there is 
    a measurable (as a function of $\omega$) sequence of paths providing the minima in the definition of $A^n_\omega(nv)$. If $\gamma^n_{A,\omega}(nv)$ is such a sequence, then
    the following relation holds for $\Pp$-a.e.\ $\omega\in\Omega$: 
    \begin{equation}
    \label{eq:gradient}
        \nabla \Lambda(v) = \lim_{n\to\infty}\frac{1}{n}\sum_{k=0}^{n-1}M(\Delta_k\gamma^n_{A,\omega}(nv)).
    \end{equation}
\end{theorem}

For a growth model interpretation, let us make an additional assumption 
\begin{equation}
    \label{eq:L-nonnegative}
    L(v)\ge 0,\quad v\in\R^d.
\end{equation}
Then we automatically have $\Lambda(v)\ge 0$.
Let us make one more assumption that this inequality is actually strict, i.e.,
\begin{equation}
\label{eq:Lambda positive}
    \Lambda(v)>0,\quad v\in\R^d.
\end{equation}
This is a nontrivial assumption. It
may be violated due to a percolation cone effect similar to that resulting in a flat edge for directed FPP limit shapes, see~\cite{Durrett-Liggett:MR0606981}. However, in multiple situations, and in particular, in our examples considered in Section~\ref{sec:Examples}, it is easy to check~\eqref{eq:Lambda positive}. 
The combination of~\eqref{eq:L-nonnegative} and~\eqref{eq:Lambda positive} will be referred to as \textbf{Assumption (L3)\namedlabel{assum:L3}{\rm (L3)}.}

For $r>0$, the random region 
\[
E_\omega(r)=\{(n,x)\in\N\times\R^d: A^n_\omega(x)\le r \} 
\]
consists of points that can be reached from the origin at a cost not exceeding~$r$. These regions grow in $r$, i.e., $0<r<r'$ implies $E_\omega(r)\subset E_\omega(r')$. Thus we obtain a random growth model. There is also an interpretation of $E_\omega(r)$ as the set of points that have been reached by a signal emitted from the origin by ``time''~$r$. One then can be interested in the asymptotic propagation of the signal and its effective front as $r\to \infty$.   

The limit shape results that describe the asymptotically linear growth of such regions are usually proved for models of First Passage Percolation (FPP) type. Although our $1+d$-dimensional model has a distinguished time coordinate, thus being of Directed FPP or Last Passage Percolation type, one can still prove a limit shape theorem:

\begin{theorem}\label{thm: limit shape}   Suppose Assumptions~\ref{assum:F},~\ref{assum:L1},~\ref{assum:L2}, 
and~\ref{assum:L3} hold. Then,  as $r\to\infty$, the set  
$\displaystyle \frac{1}{r}E(r)=\bigg\{\Big(\frac{t}{r},\frac{x}{r}\Big):\ (t,x)\in E(r)\bigg\}$ converges in the Hausdorff distance to 
the compact convex set
\begin{equation}
\label{eq:limit shape}
E_\Lambda = \overline{\big\{(t,x)\in(0,\infty)\times \R^d : t\Lambda (x/t)\le 1\big \}}.
\end{equation}
\end{theorem}
\begin{remark}
 Although the Hausdorff distance does not distinguish between a set and its closure, we prefer to define the limit shape as a closed set. 
\end{remark}

We are more interested in the properties of the limit shape than its existence, which is fairly standard.
Thus, instead of giving a rigorous argument, see \cite{Bakhtin-Dow:MR4945174}, we only recall the heuristic behind it.
According to the shape function theorem, for a large $r$, the optimal cost to reach a point $(rt,rx)$ is approximately
$rt\Lambda((rx)/(rt))=rt \Lambda(x/t)) $. Therefore, 
\[E(r)\approx \{(rt,rx):\ rt \Lambda (x/t)\le r\}=\{(rt,rx):\ t \Lambda (x/t)\le 1\}=rE_\Lambda,\]
implying the desired convergence.

\bigskip

The boundary of $E_\Lambda$ is composed of (i) points $(t,x)$ satisfying $t>0$ and $t\Lambda (x/t)=1$ and (ii) some points on the hyperplane  $H=\{(t,x):\ t=0\}$. The former may be called the effective front of signal propagation and the latter, $H_\Lambda=\partial E_\Lambda\cap H$, the footprint of the limit shape. Let us describe both parts.

Theorems~\ref{thm: main} and \ref{thm: limit shape} along with the implicit function theorem imply the following differentiability result for the effective front: 
\begin{theorem} \label{thm: limit shape C1}
Under the conditions of Theorem~\ref{thm: limit shape},  the effective front $\partial E_\Lambda \setminus H$ is a $C^1$ hypersurface diffeomorphic to $\R^d$.
\end{theorem}

For each unit vector 
$u\in\Sph^{d-1}=\{u\in\R^d:\, |u|=1 \}$, 
we define
\begin{equation}
\label{eq:linear rate lambda u}
\lambda(u)=\lim_{s\to\infty}\frac{\Lambda(su)}{s}\in (0,+\infty].
\end{equation}
Introducing 
\begin{equation}
\label{eq:G_Lambda}
G_\Lambda=\{u\in \Sph^{d-1}:\ \lambda(u)<\infty\}
\end{equation}
and using the convention $1/\infty=0$,
we can give the following description of the footprint $H_\Lambda$:
\begin{theorem} 
\label{thm:footprint} Under the conditions of Theorem~\ref{thm: limit shape},
the footprint is given by $H_\Lambda=\{0\}\times H_{\Lambda,0}$, where 
\[H_{\Lambda,0}=\Big\{ru:\ u\in \Sph^{d-1}: 0\le r \le \frac{1}{\lambda(u)}\Big\}.\]
If $G_\Lambda=\emptyset$, then $H_{\Lambda,0}=\{0\}$, and if $G_\Lambda\ne \emptyset$, then 
\[
H_{\Lambda,0}=\Big\{ru:\ u\in G_\Lambda,\ 0\le r \le \frac{1}{\lambda(u)}\Big\}.
\] 
\end{theorem}
We give the proof of this theorem in Section~\ref{sec:footprint}.

Since the set $H_\Lambda$ is convex, it has a well-defined dimension.
In Section~\ref{sec:Examples}, we consider examples where the dimension of $H_\Lambda$ equals $0$, (i.e., it is a point), $d$ (i.e., it is a $d$-dimensional ball up to a bi-Lipschitz homeomorphism) or a number strictly between those values. One example of the latter is 
provided by the Lagrangian given in~\eqref{eq:abs v1 plus squared v2}.

\section{Examples}
\label{sec:Examples}
In this section, we will discuss several examples starting with those that could be studied using the existing approaches, and emphasizing the examples where the existing results fail and one needs our new result. Our examples are accompanied with plots. Throughout this section, we will assume that $\nu=\delta_0$ in the definition of the Poisson environment. In other words, each Poisson point carries weight $\varphi=0$.

\subsection{Quadratic $L$}\label{sec:quadratic}
For all $d\in\N$,
The case of $L(v)=|v|^2$ is very similar to those studied in~\cite{BCK:MR3110798}, \cite{kickb:bakhtin2016}, \cite{Bakhtin-Li:MR3911894}. In this case, the action of a path is transformed nicely by shear transformations. 
The shear map $\Xi_v, v \in \R^d$ acting on paths $\gamma$ is given by 
\[
(\Xi_v\gamma)_k = \gamma_k + kv.
\]
The environment shear $\Xi_v^*: \Omega \to \Omega, $ is given for our model by 
\[
\Xi_v^*\omega = \{(k, x +kv,\varphi) : (k, x,\varphi) \in \omega\},\quad v\in\R^d.
\]
The stationarity in space and independence in time in the definition of the measure $\Pp$ imply that
under $\Pp$,
\begin{equation}
\label{eq:environment invariance}
\Xi_v^*\omega \overset{d}{=} \omega,\quad v\in\R^d. 
\end{equation}
In other words, $\Pp$ is preserved by the group of shear transformations.

For a path $\gamma\in\Gamma^n(0)$, we define
\begin{equation}
\label{eq:B}
B^n_\omega(v,\gamma)=A^n_{\Xi_v^*\omega}(\Xi_v\gamma).
\end{equation}
It is straightforward to compute 
\begin{align*}
B^n_\omega(v,\gamma)&=  \sum_{k=0}^{n-1} |\Delta_k( \Xi_v\gamma)|^2+\sum_{k=1}^{n-1} F_{\Xi_v^*\omega}(k,(\Xi_v\gamma)_k)
\\ &=  \sum_{k=0}^{n-1} |\Delta_k\gamma+v|^2+\sum_{k=1}^{n-1} F_{\omega}(k,\gamma_k)
\\ &=  n|v|^2+\sum_{k=0}^{n-1} |\Delta_k\gamma|^2+\sum_{k=1}^{n-1} F_{\omega}(k,\gamma_k)
\\ &= n|v|^2+A^n_{\omega}(\gamma)
\\ &= n|v|^2+B^n_{\omega}(0,\gamma).
\end{align*}
Since $\Xi_v$ establishes a bijection between $\Gamma^n(0)$ and $\Gamma^n(nv)$, 
this calculation shows that the minimizer of $A^n_{\Xi_v^*\omega}$ among paths in $\Gamma^n(nv)$ 
is the shear image of the minimizer of $A^n_{\omega}$ among paths in $\Gamma^n(0)$, and 
\begin{equation}
\label{eq:shear-inv_A}
A^n_{\Xi_v^*\omega}(nv)=n|v|^2+A^n_{\omega}(0),
\end{equation}
or, equivalently,
\begin{equation}
\label{eq:shear-inv_B}
B^n_{\omega}(v)=n|v|^2+B^n_{\omega}(0),
\end{equation}
where

\begin{equation}\label{eqn:Bnx}
    B^n_\omega(v) = \text{inf}\{B^n_\omega(v,\gamma) : \gamma \in \Gamma^n(0) \},\quad n\in \N,\ v\in\R^d.
\end{equation}
Applying the definition of the shape function and using \eqref{eq:environment invariance}, we obtain 
\begin{equation*}
    \Lambda(v)=|v|^2+C,\quad v\in\R^d,
\end{equation*}
an explicit expression up to an unknown constant $C=\Lambda(0)$.

To see that $C>0$, we can apply the following general result which also implies strict positivity of shape functions for all the other examples in this section. This allows to apply the results on limit shapes, effective fronts, and footprints from Section~\ref{sec:Main Results} to all these examples.
\begin{lemma}\label{lem:Lambda-positive}
Suppose the conditions of Theorem~\ref{thm: main}  hold and assume~\eqref{eq:L-nonnegative}. If for sufficiently small $\delta>0$,
$L_*(\delta)=\inf \{L(w): w\notin S_\delta\}>0,$
then condition~\eqref{eq:Lambda positive} holds.
\end{lemma}
We prove this lemma in Section~\ref{sec:limit-shape-proofs}.

In many examples, it is natural to define $S$ as the set of gradient discontinuities. Although $L(v)=|v|^2$ has no gradient discontinuities, we still can apply this lemma to $S=\{0\}$, $M=2v$, $H\equiv K\equiv 1$, and $u_1 = \frac{1}{2}e_1, ..., u_d = \frac{1}{2}e_d,$ $u_{d+1} = \frac{1}{d+1}(-1, ..., -1)$.
Here $e_k$ is the $k$-th vector of the standard basis in $\R^d$.

Now we can use \eqref{eq:limit shape} to find that all the points $(t,x)$ on the boundary of the limit shape 
$E_\Lambda$ that satisfy $t>0$ must also satisfy
\[
 t\Big(C+\Big|\frac{x}{t}\Big|^2\Big)= 1,
\]
or, equivalently,
\[
C\Big(t-\frac{1}{2C}\Big)^2+|x|^2=\frac{1}{4C}.
\]
Together with the origin $(0,0)$ (the only point in the footprint of the limit shape), these points form an ellipsoid centered at $(t_0,x_0)=((2C)^{-1},0)$ and invariant under spatial rotations.

\begin{figure}[ht]
    \centering
    \includegraphics[width=0.9\textwidth]{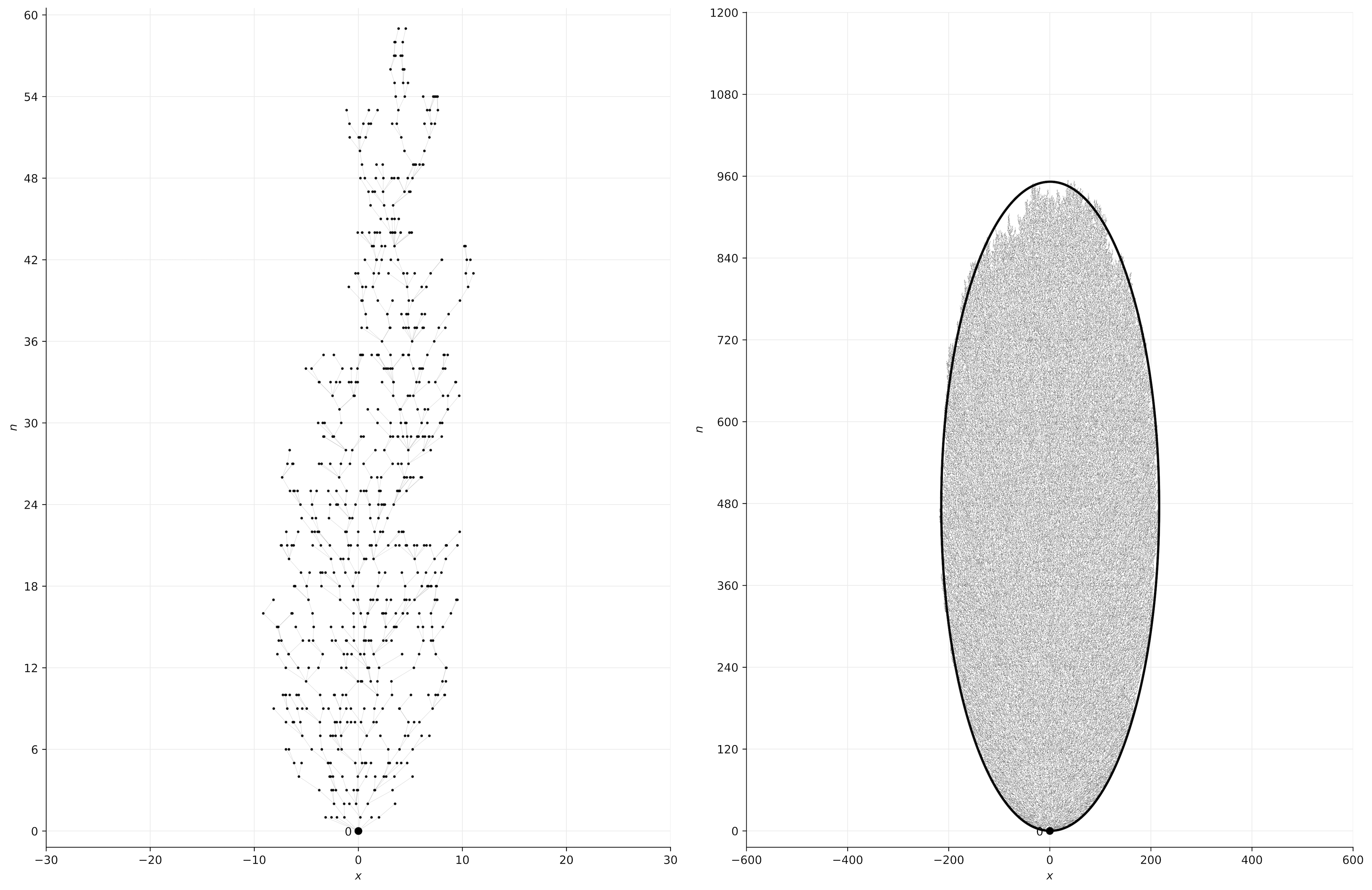} 
    \caption{Tree of minimizers and approximation to the limit shape for $L(v) = v^2$ in dimension $1+1$. On the left, points with action less than 10 are shown. On the right, points with action less than 200 are shown, along with an ellipse that is an estimate to the true limit shape.}
    \label{fig:v^2 1+1}
\end{figure}

\begin{figure}[ht]
    \centering
    \includegraphics[width=0.8\textwidth]{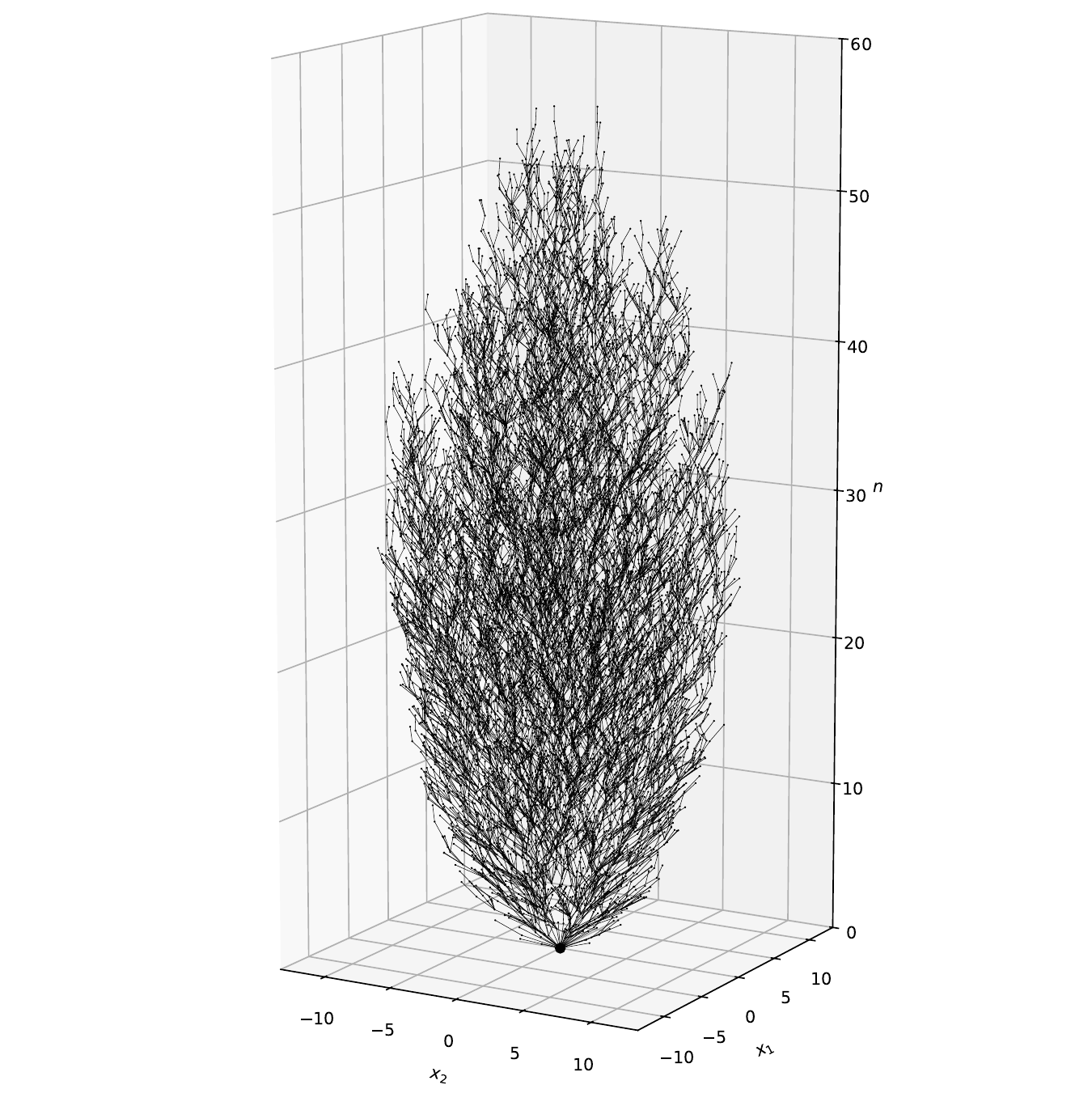} 
    \caption{Tree of minimizers and approximation to the limit shape for $L(v) = |v|^2$ in dimension $1+2$. Points with action less than 10 are shown.}
    \label{fig:v^2 1+2}
\end{figure}

Trees of minimizers and approximations to the ellipsoidal limit shape for a sampled realization of the Poisson environment are shown in Figure~\ref{fig:v^2 1+1} for $d=1$ and Figure~\ref{fig:v^2 1+2} for $d=2$. 

Obviously, the quadratic shape function is differentiable and strictly convex, and so is the boundary of the ellipsoidal limit shape.

\subsection{Polynomial $L$} If $L$ is a polynomial in $v=(v_1,\ldots, v_d)$ satisfying $\lim_{|v|\to\infty}L(v)=+\infty$
(the degree of such a polynomial must be an even positive 
number), such as  
$L(v)=|v|^4$,
then  no precise formula similar to \eqref{eq:shear-inv_A} or \eqref{eq:shear-inv_B} exists (except the quadratic case). However, as~\cite{BakhtinDow-AIHP1529} shows, the shear-transformed action~$B^n$ defined by \eqref{eq:B} is still useful.
To see the advantage of working with $B^n$, we note that
$\Gamma^n(nv)$, the set of admissible paths one minimizes over to obtain $A^n(nv)$,  depends on $v$ and thus explore different domains of the random environment. Meanwhile in the definition of $B^n_\omega(v)$,
no matter what $v$ is, the minimization is taken over the same collection of paths $\Gamma^n(0)$ in the same realization of the random environment.

\begin{figure}[ht]
    \centering
    \includegraphics[width=0.8\textwidth]{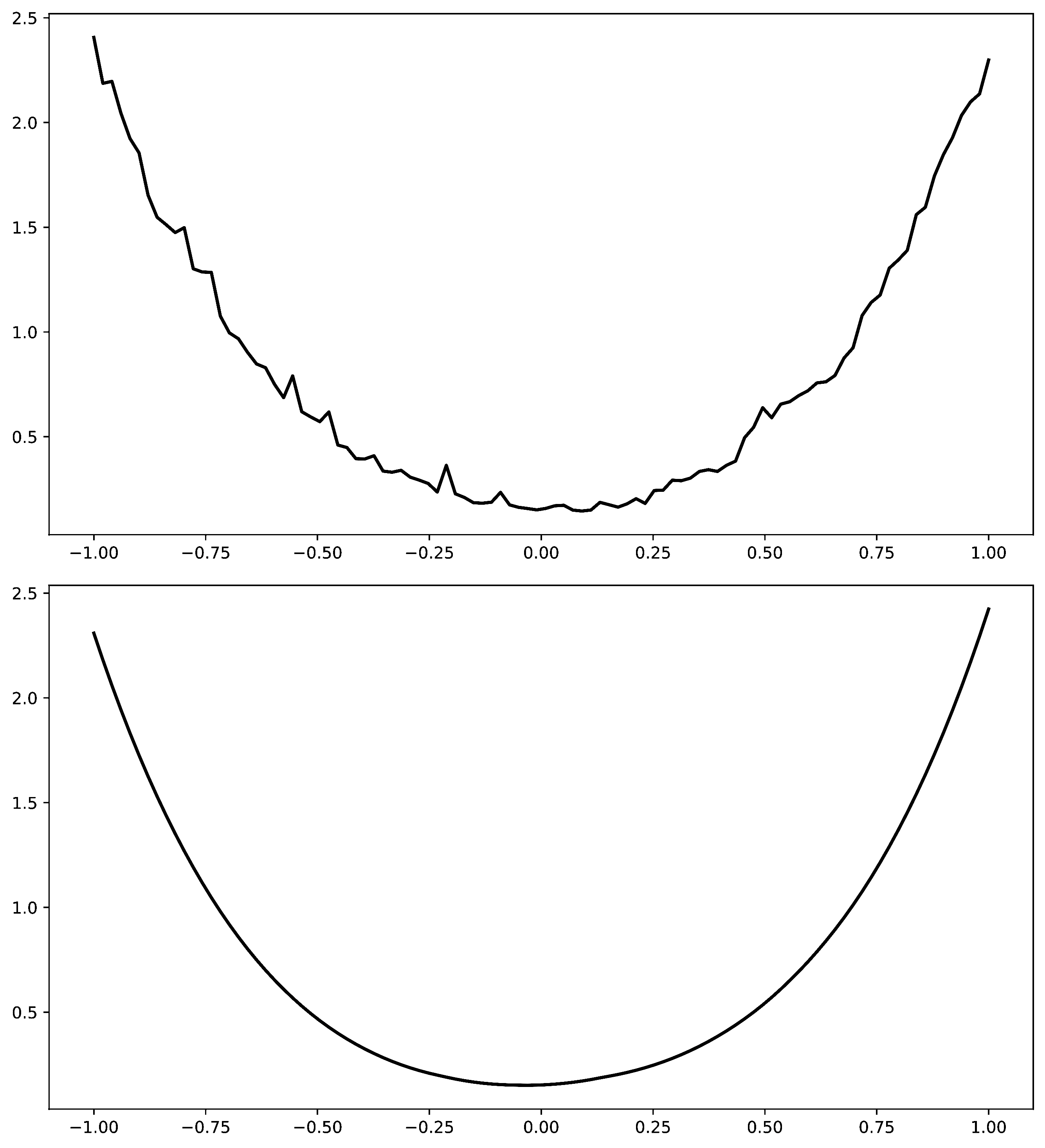} 
    \caption{Approximation to the shape function for $L(v) = v^4$ in dimension $1+1$. Above is $\frac{A^n_\omega(nv)}{n}$, while below is $\frac{B^n_\omega(v)}{n}$, for $n = 200$. }
    \label{fig:v^4}
\end{figure}

In our model, the collection of paths is a discrete set, so if a path provides the minimum in the definition of $B^n_\omega(v)$  for some $v$, it also serves as a minimizer for nearby values of $v$. Thus, $B^n_\omega$ is piecewise polynomial. For example,
if $d=1$ and $L(v)=v^4$, the graph of $B^n_\omega$ can be described 
as a concatenation of several convex quartics, each one corresponding to a path optimal for a whole range of values of~$v$, see Figure~\ref{fig:v^4}.
In the case of $L(v)=v^2$, the entire graph would consist of just one parabola since the minimizer realizing $B^n_\omega(0)$ is also the minimizer realizing $B^n_\omega(v)=nv^2+B^n_\omega(0)$. 

The method of~\cite{BakhtinDow-AIHP1529} based on the analysis of $B^n$ can be used verbatim to prove the differentiability of the shape function $\Lambda$ in this case. Moreover, 
\begin{equation}
\label{eq:gradient in C2 case}
\nabla\Lambda(v)\aseq\lim_{n\to\infty}\frac{1}{n}\sum_{k=0}^{n-1}\nabla L(\Delta_k\gamma^n(v)),
\end{equation}
where $\gamma^n(v)$ is a path that provides the minimum in the definition of $A^n(nv)$ in~\eqref{eqn:Anx}.
Note that if $L(v)=|v|^2$, this formula can be checked directly since 
\[
\frac{1}{n}\sum_{k=0}^{n-1}\nabla L(\Delta_k\gamma^n(v))= \frac{1}{n}\sum_{k=0}^{n-1}(2\Delta_k\gamma^n(v))=2v=\nabla\Lambda(v).
\]
\begin{figure}[ht]
    \centering
    \includegraphics[width=0.9\textwidth]{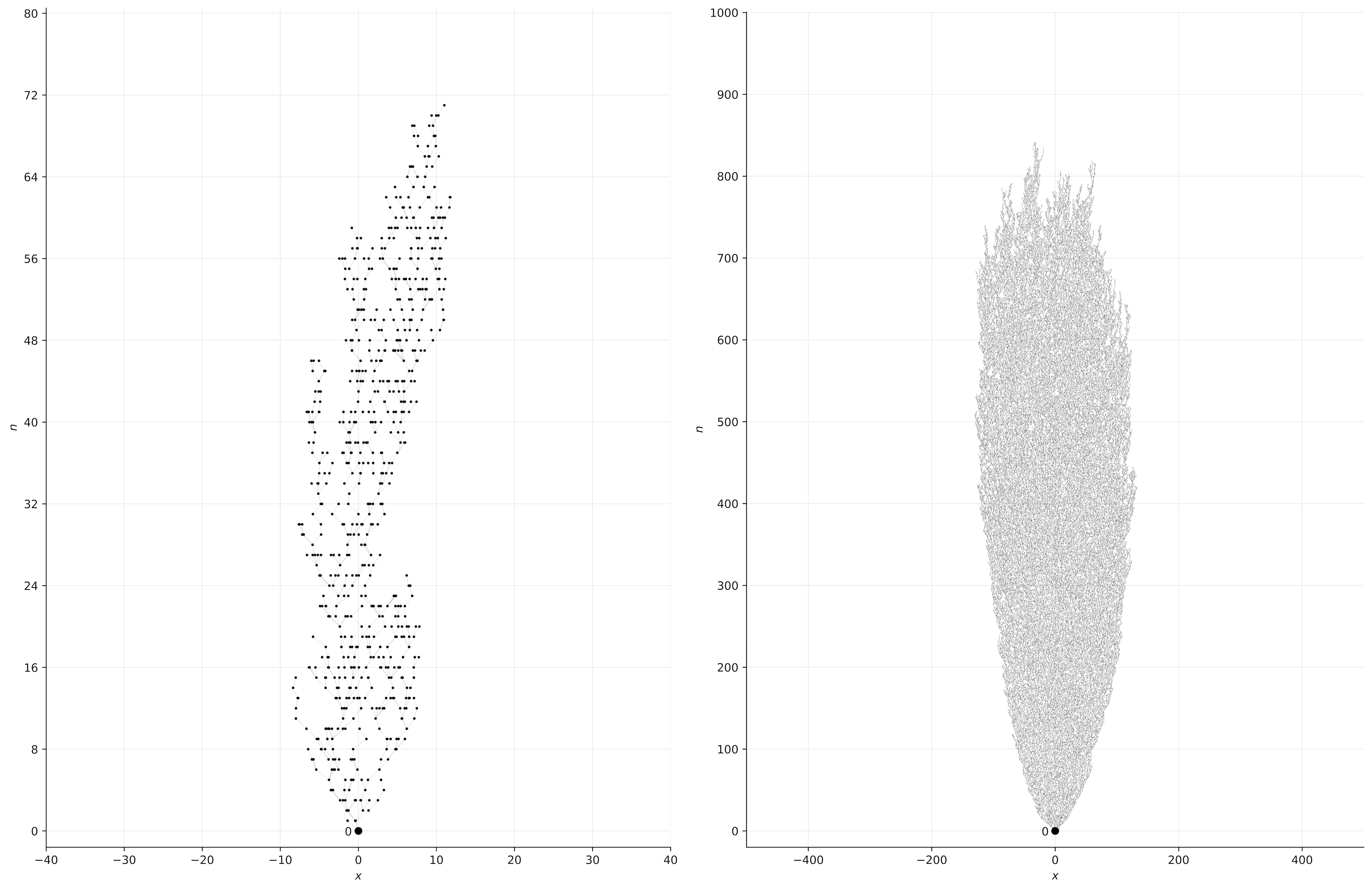}
    \caption{Tree of minimizers and approximation to the limit shape for $L(v) = v^4$ in dimension $1+1$. On the left, points with action less than 10 are shown. On the right, points with action less than 100 are shown.}
    \label{fig:v^4 1+1}
\end{figure}

It can easily be shown that for polynomial $L$ considered in this section, the growth of $\Lambda$ at infinity is superlinear. This implies that the footprint of the limit shape consists of a single point at the origin.

\subsection{$L$ is a norm} If $L$ is a norm such as the Euclidean norm, then $\nabla L$ has a discontinuity at $0$ and there could be entire hypersurfaces of discontinuity, e.g., for $L^1$ and $L^\infty$ norms, so the methods from the existing literature do not apply directly but the results of the present paper do apply.

If $d=1$ and $L(v)=|v|$, then the conditions of our Theorem~\ref{thm: main} apply if we choose $S=\{0\}$, $u_1=1$, $u_2=-1$, $H\equiv 0$, $\alpha=1$, $r_0 = 1$, $\epsilon_0 = 1$, $K(v)=2$, $c_0=2$, $c_1=0$,  $\eta_R\equiv1$, $k_R\equiv 2$, 
\begin{equation*}
M(v)=\sgn v=\begin{cases}
    1,& v>0,\\ 0,& v=0\\ -1,& v<0.
\end{cases}    
\end{equation*}
We note that it is natural here and in other applications to choose $M$ to be the weak derivative or gradient of $L$. If $L$ is piecewise differentiable, then this means defining $M$ to be the classical gradient where it is well-defined and setting it to $0$ at gradient discontinuities. Applying Theorem~\ref{thm: main}, we obtain that the shape function $\Lambda$ belongs to $C^1(\R)$ and, moreover, formula \eqref{eq:gradient} which is a generalization of \eqref{eq:gradient in C2 case} takes a particularly simple form:
\begin{equation}
\label{eq:gradient-for-abs-L}
\Lambda'(v)\aseq\lim_{n\to\infty}\frac{1}{n}\sum_{k=0}^{n-1} \sgn(\Delta_k\gamma^n(v)). 
\end{equation}
In words, the derivative of $\Lambda(v)$ is computed through the balance of the number of steps that the optimal path makes to the right and to the left. The simplicity of the formula gives hope that it will be used to prove strict convexity of $\Lambda$ near~$0$. The behavior is different outside of $[-1,1]$:
\begin{theorem}
\label{thm:flat edges}
If $d=1$, $L(v)=|v|$, and $\nu=\delta_0$, then
for all $|v|\ge 1$, we have $\Lambda(v)=|v|$.
\end{theorem}
\begin{proof} The triangle inequality implies that for all $v$, all $\gamma\in\Gamma^n(nv)$, and all $\omega\in\Omega$, $A^n_\omega(\gamma)\ge n|v|$, implying $\Lambda(v)\ge |v|$. 

Now suppose $v>1$. For each $n\in\N$, we can define the following path $\gamma^n$ recursively:
$\gamma_0=0$, 
\[\gamma_k=\inf\{x> \gamma_{k-1}: (k,x)\in\pi\omega\},\quad k=1,\ldots, n-1,\]
and $\gamma_n=nv$. Then $(\gamma_{k}-\gamma_{k-1})_{k=1,\ldots,n-1}$ is an i.i.d.~collection of rate $1$ exponential r.v.'s.
Applying the strong law of large numbers, we obtain that with probability $1$, for sufficiently large $n$, $\gamma_{n-1}<nv$,
so that $A^n(\gamma)=nv$ implying $\Lambda(v)\le v$. By continuity, this inequality extends to $v=1$ and by symmetry to $v\le -1$. Combining this with the lower estimate obtained for all $v$, we complete the proof.    
\end{proof}

The strictly linear behavior of $\Lambda$ outside of $[-1,1]$ results in an interesting limit shape for the growth model. 
\begin{theorem} Under the conditions of Theorem~\ref{thm:flat edges}, the effective front is $C^1$. It has two flat pieces: $(0,1]\times\{1\} $ and $(0,1]\times\{-1\}$.  
\end{theorem}
\begin{proof}
Since Lemma~\ref{lem:Lambda-positive} with $L_*(\delta)=\delta$ implies~\eqref{eq:Lambda positive}, we can apply Theorem~\ref{thm: limit shape C1} to derive that the front is $C^1$. 
Moreover, 
Theorem~\ref{thm: limit shape}
gives an equation for the effective front, and 
from Theorem~\ref{thm:flat edges} we obtain that points $(t,x)$ satisfying $|x|= 1$ and $t\in(0,1]$, satisfy this equation: $t\Lambda(x/t) = t|x|/t= 1$.
\end{proof}

We conjecture that the central part of the front complementary to the union of these two flat edges is a strictly convex curve.
\begin{figure}[htbp]
    \centering
    \includegraphics[width=0.9\textwidth]{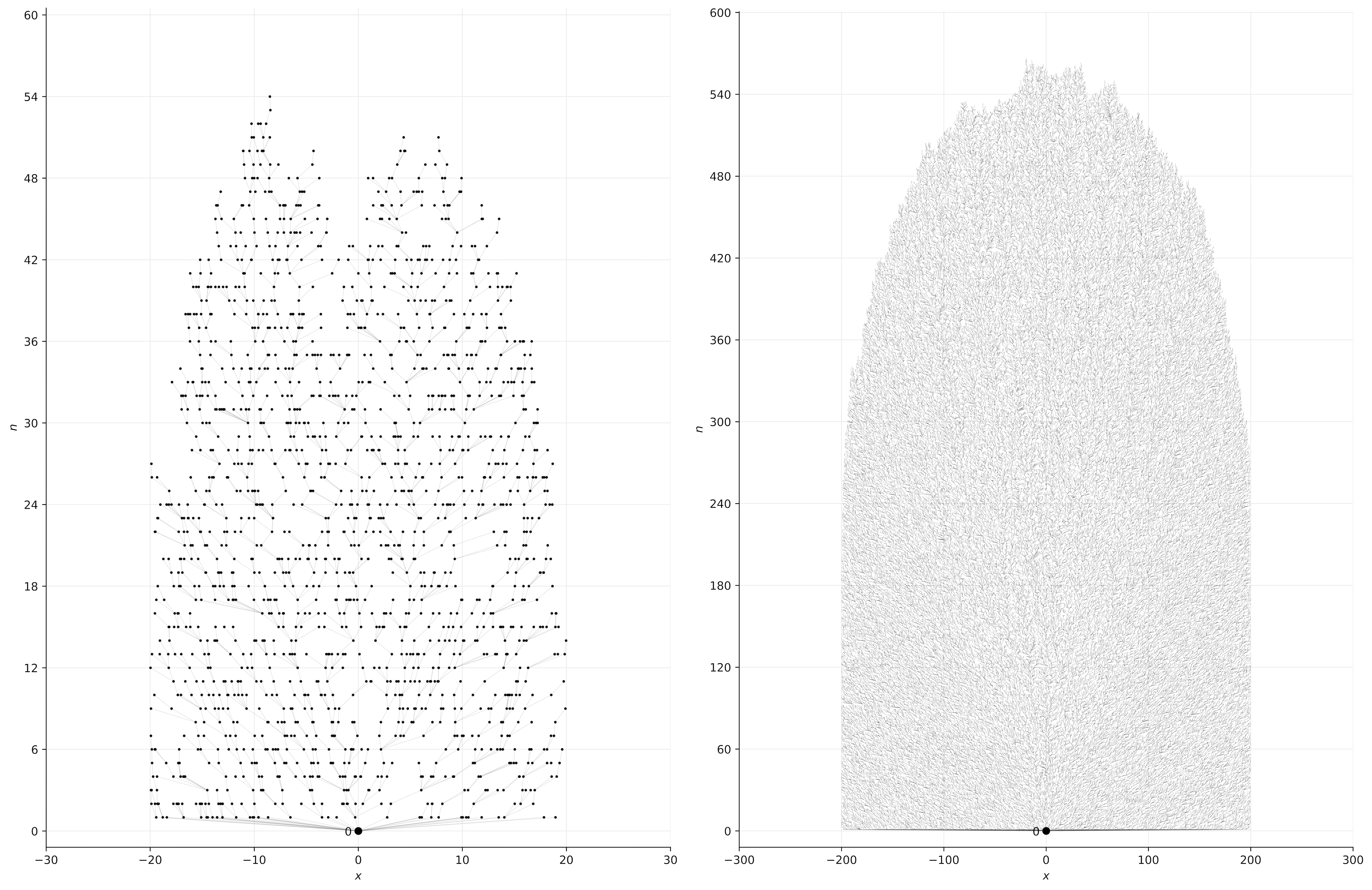}
    \caption{Tree of minimizers and approximation to the limit shape for $L(v) = |v|$ in dimension $1+1$. On the left, points with action less than 20 are shown. On the right, points with action less than 200 are shown.}
    \label{fig:|v| 1+1}
\end{figure}

If $d \ge 2 $ and $L(v) = |v|$, the conditions of 
Theorem~\ref{thm: main} apply with $S = \{0\},$ $u_1 = \frac{1}{2}e_1, ..., u_d = \frac{1}{2}e_d,$ $u_{d+1} = \frac{1}{d+1}(-1, ..., -1),$ $M(v) = \frac{v}{|v|}\ONE_{v\neq 0},$ $K \equiv 2,$ $H(v) = \frac{2}{|v|}\ONE_{v\neq 0},$ $\alpha = 1,$ $\theta = 0,$ $\eps_0 = 1,$ $r_0 = 1,$ $c_0 = 2,$ $c_1 = 1,$ $c_2 = 2,$  $\eta_R = d,$   $k_R = m\big(B_d(0, 1)\big)$. 
Since  the distribution of the environment is invariant under rotation and the Lagrangian is radially symmetric, the limit shape will also be radially symmetric. 
Define the generating curve as the set 
\[
\{\big(0, (1, 0, ..., 0)\big)\} \cup \{\big(t, (x, 0, ..., 0) \big) : x  \geq 0, t \Lambda\big((x/t, 0, ..., 0)\big) = 1\}. 
\]
\begin{theorem}\label{thm: d>1}
    Suppose $d\ge 2$. If $ L(v) = |v| $ and $\nu = \delta_0$, the footprint of the limit shape is the closed $d$-dimensional unit ball. Furthermore, the closure of the effective front consists of points swept out by spatial rotations of the generating curve, and the generating curve is not linear at $t= 0$. 
\end{theorem}
We give a proof of this theorem in Section~\ref{sec:limit-shape-proofs}.

\begin{figure}[ht]
    \centering
    \includegraphics[width=0.8\textwidth]{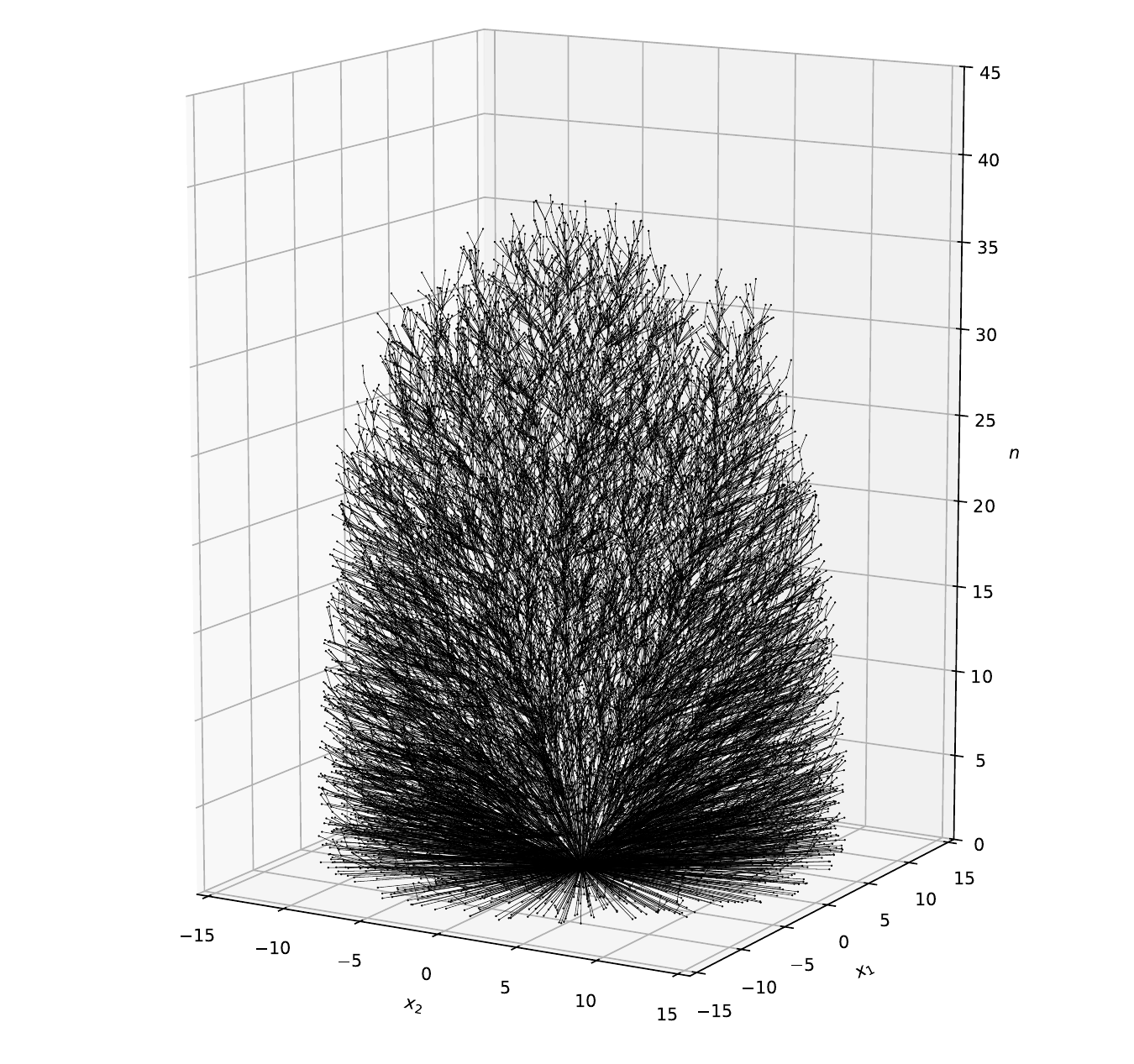} 
    \caption{Tree of minimizers and approximation to the limit shape for $L(v) = |v|$ in dimension $1+2$. Points with action less than 15 are shown.}
    \label{fig:|v| 1+2}
\end{figure}

Comparing limit shapes for the cases we have considered, we notice that their footprints look different.
In the polynomial $L$ case, the footprint consists of a single point, the origin. In the case of $L(v)=|v|$ in $d=1$, the footprint is the segment $\{0\}\times[-1,1]$. These features are in agreement with  Theorem~\ref{thm:footprint}.
Let us consider an example combining those features.

\subsection{Mixed $L$}
 Let us take arbitrary dimensions $d\in \N$ and $m\ge 0$  satisfying $0\le m\le d$ and define
\begin{align}
\label{eq:mixed L}
L(v)&=|v_1|+\ldots+|v_{m}| + v_{m+1}^2+\ldots+v_d^2
\\ \notag &= |\Pi_m v|_1+|\Pi_m^\perp v|^2,
\end{align}
where $\Pi_m v =(v_1,\ldots,v_m,0,\ldots,0)$ is the orthogonal projection on the vector space $H_m$ spanned by the first~$m$~coordinates, $\Pi_m^\perp v =v-\Pi_m v$ is the orthogonal projection on the complementary coordinates, $|\cdot|_1$ is the $L^1$ norm, and $|\cdot|$ is the Euclidean (i.e., $L^2$) norm. We allow for the degenerate cases $m=0$ or $m=d$. 

For each $i=1,\ldots,m$, each point on the hyperplane $\{x: x_i=0\}$ is a gradient discontinuity, but it is easy to see that Theorem~\ref{thm: limit shape C1} applies and 
 implies the following: 
\begin{theorem} 
\label{thm:front for mixed}
For any $d$ and $m$,
the effective front for the Poissonian growth model defined by $L$ given in \eqref{eq:mixed L} is a $C^1$ hypersurface.
\end{theorem}
In fact, using shears similarly to Section~\ref{sec:quadratic} but only in the direction of the last $d-m$ coordinates, we obtain
\begin{equation}
\label{eq:Lambda for mixed}
\Lambda(v)=\Lambda(\Pi_m v)+|\Pi_m^\perp v|^2.    
\end{equation}
We claim that there is a constant $C>0$ such that
\begin{equation}
\label{eq:mixed Lambda linear growth}
  \Lambda(v)\le C+ |v|_1,\quad v\in H_m. 
\end{equation}
Relations \eqref{eq:Lambda for mixed} and \eqref{eq:mixed Lambda linear growth} imply that 
the set  $G_\Lambda$ defined in
\eqref{eq:G_Lambda} equals $\Sph^{d-1}\cap H_m$. Thus Theorem~\ref{thm:footprint} implies the following:
\begin{theorem} \label{thm:footprint for mixed}
For any $d$ and $m$,
the footprint of the limit shape for the Poissonian growth model defined by $L$ given in \eqref{eq:mixed L} is an $m$-dimensional convex set.  
\end{theorem}

It remains to establish~\eqref{eq:mixed Lambda linear growth} for $m\ge 1$.
To that end, we first fix $v\in H_m$ and construct a path $\gamma\in\Gamma^n(nv)$  providing an upper estimate on $A^n(nv)$ for sufficiently large $n$.

We must have $\gamma_0=0$.
Let's choose $\gamma_1,\gamma_2,\ldots, \gamma_{n-1}$ recursively. If $\gamma_{k-1}$ has been chosen, then we set $\gamma_k$ to be the point $x$ minimizing
$|\Pi_m(x-\gamma_{k-1})|_1$ among the Poissonian points $(k,x)$ satisfying 
$|\Pi_m^\perp x|\le 1/2$. Then $|\Pi_m\Delta_k\gamma|_1$ is a
finite mean i.i.d.\  sequence
and $|\Pi_m^\perp\Delta_k\gamma|^2\le 1$ for all $k$. 
The strong law of large numbers implies that there is a constant $C'$ (independent of $v$) such that with probability $1$, for sufficiently large $n$, 
\begin{equation}
\label{eq:sum Pi_m}
\sum_{k=0}^{n-2} |\Pi_m\Delta_k\gamma|_1\le C'n.
\end{equation}
This implies $|\Pi_m\gamma_{n-1}|_1\le C'n$. 
Thus, setting $\gamma_n=nv$,  we obtain
$|\Pi_m \Delta_{n-1}\gamma|_1\le (C'+|v|_1)n$, which implies
\begin{equation}
\sum_{k=0}^{n-1} |\Pi_m\Delta_k\gamma|_1\le (2C'+|v|_1)n.
\end{equation}
Since 
\[
\sum_{k=0}^{n-1} |\Pi_m^\perp\Delta_k\gamma|^2\le n,
\]
we obtain~\eqref{eq:mixed Lambda linear growth} with $C=2C'+1$.

\begin{figure}[ht]
    \centering
    \includegraphics[width=\textwidth]{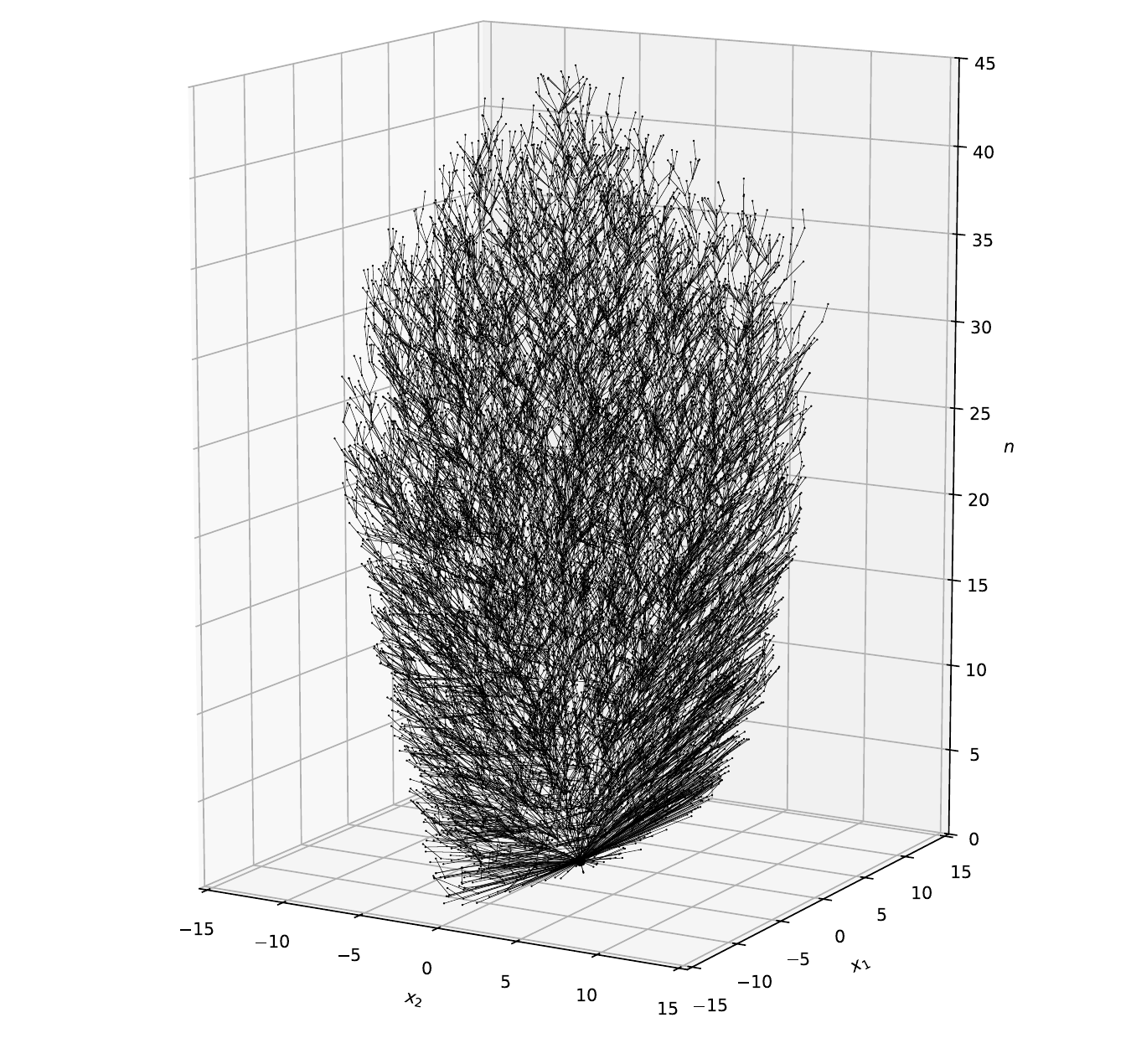} 
    \caption{Tree of minimizers and approximation to the limit shape for $L(v) = |v_1| + v_2^2$ in dimension $1+2$. Points with action less than 15 are shown.}
    \label{fig:mixed L}
\end{figure}

Theorems~\ref{thm:front for mixed} and \ref{thm:footprint for mixed} apply directly (with $d=2$ and $m=1$) to the Lagrangian defined in~\eqref{eq:abs v1 plus squared v2}. The resulting effective front is $C^1$ and the footprint is a $1$-dimensional convex shape, i.e., a line segment, see Figure ~\ref{fig:mixed L}.

\section{Proof of the Main Result}\label{sec:main_proofs}

\subsection{The Sheared Environment} 
The family of shear transformations $\Xi_v$ and $\Xi^*_v$, the transformed action $B^n_\omega(v,\cdot)$ and transformed minimal action $B^n_\omega(v)$  were introduced in Section~\ref{sec:quadratic}. Although in general there is no simple algebraic expression for the transformed action similar to~\eqref{eq:shear-inv_B} for the quadratic case, it will still be useful to work in the sheared environment.

\begin{proposition} \label{prop:min}
For all realizations $\omega$ of the environment, for all $ n \in \mathbb{N}$, for all $v \in \mathbb{R^d}$, there exists a minimizer $\gamma_A^n(nv)$ providing the minimum in the definition~\eqref{eqn:Anx} with $x=nv$ 
and a minimizer $\gamma_B^n(v)$ for the functional $B^n_\omega(v)$ defined by~\eqref{eqn:Bnx}, these minimizers can be chosen measurably, and
\begin{equation}\label{eq:dist}
    (B^n_{\omega}(v), \Xi_v\gamma^n_B(v))_{n\in\N}  \overset{d}{=} (A^n_{\omega}(nv), \gamma^n_A(nv))_{n\in\N}.
\end{equation}
\begin{proof}
The claim will follow from the distributional invariance of the Poisson point process under shear transformations \eqref{eq:environment invariance} and the definition of the transformed optimal action $B^n_\omega(v)$ in
\eqref{eqn:Bnx} once we choose the minimizers measurably.

Since $L(v) \to \infty$ as $|v| \to \infty$ and the Poisson process is locally finite, we only need to consider finitely many paths as possible minimizers, and so the minima are well defined.

We can view paths $\gamma \in \Gamma^n(x)$ as elements of $\R^{d \times (n+1)}$. Then we can order the (finitely many) minimizing paths using lexicographical order. Choosing the first path in the order defines a measurable selection of the minimizers.
\end{proof}
\end{proposition}

\begin{remark}\label{rem:Lambda-via-B}
Proposition~\ref{prop:min} implies almost surely, 
\begin{equation*}
    \Lambda(v) = \lim_{n \to \infty} \frac{1}{n}B^n_{\omega}(v). 
\end{equation*}
This representation is the most important ingredient in the proof of  Theorem~\ref{thm: main}. 
\end{remark}

\begin{remark} Our approach allows one to derive differentiability  of $\Lambda$ by working with more general transformations of the form
\[
(\Xi_v \gamma)_k = \gamma_k+nf(k/n)v, 
\]
where $f$, instead of the identity map, can be taken to be any nondecreasing Lipschitz function $f:[0,1]\to[0,1]$ satisfying $f(0)=0$, $f(1)=1$. 
Moreover, different choices of $f$ result in different versions of formula~\eqref{eq:gradient}. This results in a certain regularity of the behavior of  $M$ along the optimal paths, but we do not pursue this direction in this paper.
\end{remark}

\subsection{Approximate linear domination  lemma}
Another important ingredient in our proof of differentiability of $\Lambda$ is the following lemma, which is a slight generalization of Lemma~3.3 in \cite{bakhtin2023differentiability}. It derives differentiability of a convex function at a point from an approximate linear dominance bound \eqref{limit}---\eqref{eq:lemma33-assumption}.
\begin{lemma}\label{lem: conv}

Let $O\subset\mathbb R^d$ be an open ball centered at $x_0$, and $f$ a convex function on $O$. Let $(f_n)_{n\in\mathbb N}$ be a sequence of functions from $O$ to $\mathbb R$. Suppose that there exist sequences $(x_{m}^{(i)})_{m\in\mathbb N}$, $i = 1, ..., d+1$, where $x_m^{(i)} \in O\setminus\{x_0\}$, non-negative sequences $(h_{n,m}^{(i)})_{n, m \in \N}$ in $\R$, a sequence $(\xi_n)_{n \in \N}$ in $\R^d$, and $d+1$ non-zero vectors $u_1, \dots, u_{d+1} \in \mathbb{R}^d$ that satisfy the following conditions.
\begin{enumerate}
\item For all $i \in \{1, ..., d+1\}$, 
\begin{equation}
\label{eq:x_m-to-x_0}
    \lim_{m \to \infty }x_m^{(i)} = x_0.
\end{equation}
\item  For all $i\in \{1, ..., d+1\}$,
\begin{equation}
\label{eq:norm_x_m-conv-to-u}
\lim_{m\to\infty} \frac{x_{m}^{(i)}-x_0}{|x_{m}^{(i)}-x_0|} = \frac{u_i}{|u_i|}.
\end{equation}
\item \label{it:0-in-conv-hull}
The interior of the convex hull of $u_1, ..., u_{d+1}$ contains $0$.

\item For all $m \in \N,$ and $i \in \{1, ..., d+1\}$,
\begin{equation}
\label{eq:f_n-conv-to-f-at-x_m}
    \lim_{n \to \infty}f_n(x_m^{(i)}) = f(x_m^{(i)})
\end{equation}
and 
\begin{equation}
\label{eq:f_n-conv-to-f-at-x_0}
    \lim_{n \to \infty}f_n(x_0) = f(x_0).
\end{equation}
\item For all $i \in \{1, ..., d+1\}$,
\begin{equation}\label{limit}
    \lim_{m\to\infty}\limsup_{n \to \infty}\frac{h_{n, m}^{(i)}}{|x_{m}^{(i)}-x_0|}=0.
\end{equation}

\item For all $m, n \in \mathbb{N}$ and $i \in \{1, ..., d+1\}$,
\begin{equation}\label{eq:lemma33-assumption}
f_n(x_{m}^{(i)})-f_n(x_0) \le \langle \xi_n, x_{m}^{(i)}-x_0\rangle+h_{n, m}^{(i)}.
\end{equation}

\end{enumerate}

Then $f$ is differentiable at $x_0$, the sequence $(\xi_n)_{n\in\mathbb N}$ converges, and
\begin{equation}\label{eq:lemma33-conclusion}
\nabla f(x_0)=\lim_{n\to\infty}\xi_n.
\end{equation}
\end{lemma}

\begin{proof}
We recall the definition of the subdifferential of $f$:
\[
\partial^\vee f(x)=\bigl\{\xi\in\mathbb{R}^d : f(y)-f(x)\ge\langle \xi,y-x\rangle \text{ for all } y\in O\bigr\}.
\]
By Theorem 25.1 in \cite{Rockafellar1970}, the set $\partial^\vee f(x)$ is convex and non-empty, and  $f$ is differentiable at $x$ if and only if $\partial^\vee f(x)$ contains exactly one point. This is the condition we will check.

Let's prove that $|\xi_n|$ is bounded. If not, then there is a subsequence $\xi_{n'}$ and $e\in\R^d\setminus\{0\}$ satisfying  $|\xi_{n'}|\to \infty$ and $\xi_{n'}/|\xi_{n'}|\to e$. Condition~\ref{it:0-in-conv-hull} implies that there is $i$ such that $\langle e,u_i\rangle <0$. Then there is a constant $c_1>0$ such that 
$\langle\xi_{n'}/|\xi_{n'}|, u_i\rangle<-c_1$ for large $n'$, which in turn implies, due to~\eqref{eq:x_m-to-x_0} and~\eqref{eq:norm_x_m-conv-to-u}  that there are a sufficiently large $m\in\N$ and a sufficiently small $c_2>0$ such that $\langle\xi_{n'}/|\xi_{n'}|, x_{m}^{(i)}-x_0\rangle<-c_2$ for large $n'$.
By~\eqref{eq:f_n-conv-to-f-at-x_m}---\eqref{eq:lemma33-assumption}, we
obtain that $|\xi_{n'}|\langle\xi_{n'}/|\xi_{n'}|, x_{m}^{(i)}-x_0\rangle$ is bounded, and thus so is $|\xi_{n'}|$.

Let $G\subset \mathbb R^d$ 
denote the set of limit points of the sequence $(\xi_n)_{n\in\mathbb N}$, and let $\xi^*\in G$ be a limit point of a subsequence $(\xi_{n_k})_{k\in\mathbb N}$.

Fix an index $i \in \{1, \dots, d+1\}$ and the corresponding sequence $(x_{m}^{(i)})_{m\in\mathbb N}$. Taking $k\to\infty$ in
\[
f_{n_k}(x_{m}^{(i)})-f_{n_k}(x_0) \le \langle \xi_{n_k},x_{m}^{(i)}-x_0\rangle+h_{n_k,m}^{(i)},
\]
we obtain
\begin{equation}\label{eq:lemma33-limit}
f(x_{m}^{(i)})-f(x_0) \le \langle \xi^*,x_{m}^{(i)}-x_0\rangle+\limsup_{n\to\infty}h_{n,m}^{(i)}.
\end{equation}

Let $p\in\partial^\vee f(x_0)$. Then
\[
f(x_{m}^{(i)})-f(x_0)\ge \langle p,x_{m}^{(i)}-x_0\rangle,
\]
and together with \eqref{eq:lemma33-limit} this yields
\begin{equation}\label{eq:lemma33-key}
0 \le \langle \xi^*-p, x_{m}^{(i)}-x_0\rangle+\limsup_{n\to\infty}h_{n, m}^{(i)}.
\end{equation}

Dividing \eqref{eq:lemma33-key} by $|x_{m}^{(i)}-x_0|$ and letting $m\to\infty$, using \eqref{limit} gives
\[
0\le \left\langle \xi^*-p, \frac{u_i}{|u_i|} \right\rangle,
\]
which implies
\[
0\le \langle \xi^*-p, u_i \rangle,\quad i =1, \dots, d+1.
\]

Due to Assumption~\ref{it:0-in-conv-hull}, there are  constants $\alpha_1, \dots, \alpha_{d+1} > 0$ such that
\[
\sum_{i=1}^{d+1} \alpha_i u_i = 0.
\]

Taking the inner product of $\xi^*-p$ with this linear combination yields
\[
0 = \langle \xi^*-p, 0 \rangle = \sum_{i=1}^{d+1} \alpha_i \langle \xi^*-p, u_i \rangle.
\]

Since each term $\alpha_i \langle \xi^*-p, u_i \rangle$ is nonnegative, the sum can only be zero if
\[
\langle \xi^*-p, u_i \rangle = 0 \quad \text{for all } i \in \{1, \dots, d+1\}.
\]

Because $0$ is an interior point of the convex hull of $u_1, \dots, u_{d+1}$, these $d+1$ vectors span $\mathbb{R}^d$. Since $\xi^*-p$ is orthogonal to all $u_i$, we have $\xi^*-p=0$, i.e., 
$\xi^* = p.$

Since $\xi^*\in G$ and $p\in\partial^\vee f(x_0)$ were arbitrary, we  conclude that
$G=\partial^\vee f(x_0)$
consists of a single point, so $\lim_{n\to\infty}\xi_n$ is well-defined and
\[
\partial^\vee f(x_0) = \left\{\lim_{n\to\infty}\xi_n\right\}.
\]
Therefore, $f$ is differentiable at $x_0$ and \eqref{eq:lemma33-conclusion} holds.

\end{proof}

\begin{remark}
    In practice, we will often define $h_{n, m}^{(i)} = h_n(x_m^{(i)})$ for some sequence of functions $h_n : O \to \R^+$. 
\end{remark}

\subsection{Proof of Theorem~\ref{thm: main}}

Our first lemma does not involve randomness.
To state it, for every subset of the unit sphere $\Sph^{d-1}$, we define its aperture as the maximum angle between its elements.

\begin{lemma}\label{Yao}
Let $V_0, V_1, \dots, V_N$ be locally finite non-empty subsets of $\mathbb{R}^d$. For a path $\gamma = (\gamma_0, \gamma_1, ..., \gamma_N)$ with $\gamma_i \in V_i$, define the special action
\begin{equation}\label{cost}
    \bar{A}(\gamma) = \sum_{i=0}^{N-1} |\gamma_{i+1} - \gamma_i|,
\end{equation}
and let $\gamma^* = (\gamma^*_0, \gamma^*_1, ..., \gamma^*_N)$ minimize \eqref{cost}. 

Suppose the unit sphere $\Sph^{d-1}$ is covered by a collection of $m$ spherical caps $\{K_1, \dots, K_m\}$ with aperture not exceeding some $\theta < \frac{\pi}{3}$.
Let $\bar{\Gamma}$ be the set of all valid restricted paths $\gamma = (\gamma_1, \dots, \gamma_N)$ starting at $\gamma_0 = \gamma^*_0$ such that for each $i \in \{0, \dots, N-1\}$, the next point $\gamma_{i+1}$ minimizes $|\gamma_{i+1} - \gamma_i|$ among points in $V_{i+1}$ along the direction of \emph{some} cap $K_l$. Here we adopt the convention that $0$ is along the direction of all caps. Then, the minimum special action among all valid restricted paths satisfies:
\[
\min_{\gamma \in \bar{\Gamma}} \bar{A}(\gamma) \le c_{\theta} \bar{A}(\gamma^*),
\]where
\begin{equation}\label{theta}
c_{\theta} = \frac{1}{1 - 2\sin(\theta/2)} > 0
\end{equation}
\end{lemma}

\begin{proof}
We construct a specific restricted path $\gamma \in \bar\Gamma$ sequentially. We start with $\gamma_0 = \gamma^*_0$.
Then suppose the path has reached a point $\gamma_i \in V_i$. Let $u_i = \gamma^*_{i+1} - \gamma_i$ be the vector pointing to the next point $\gamma^*_{i+1}$ in the true optimal path, and let $K_l$ be a spherical cap containing the direction vector $\frac{u_i}{|u_i|}$. Among points $z \in V_{i+1}$ such that $\frac{z-\gamma_i}{|z-\gamma_i|} \in K_l$, we choose $\gamma_{i+1}$ to minimize $|z - \gamma_i|$. Note that this choice is well-defined, since $\frac{\gamma^*_{i+1}-\gamma_i}{|\gamma^*_{i+1}-\gamma_i|} \in K_l$.

Let $v_i = \gamma_{i+1} - \gamma_i$ denote this chosen restricted step. By the definition of our choice, it satisfies:
\begin{enumerate}
    \item $|v_i| \le |u_i|$
    \item $\angle(u_i, v_i) \le \theta$
\end{enumerate}

Then we have, using the law of cosines, 
\begin{align*}
|u_i - v_i|^2 
&\le |u_i|^2 + |v_i|^2 - 2|u_i||v_i|\cos\theta \\
&= |u_i|^2 + |v_i|^2 - 2|u_i||v_i|\left(1 - 2\sin^2\left(\frac{\theta}{2}\right)\right) \\
&= \left(|u_i| - |v_i|\right)^2 + 4\left(|u_i| - |v_i|\right)|v_i|\sin^2\left(\frac{\theta}{2}\right) + 4|v_i|^2\sin^2\left(\frac{\theta}{2}\right) \\
&\le \left(|u_i| - |v_i|\right)^2 + 4\left(|u_i| - |v_i|\right)|v_i|\sin\left(\frac{\theta}{2}\right) + 4|v_i|^2\sin^2\left(\frac{\theta}{2}\right) \\
&= \left[ \left(|u_i| - |v_i|\right) + 2|v_i|\sin\left(\frac{\theta}{2}\right) \right]^2
\end{align*}
Taking square roots gives
\begin{equation}\label{eq:yao_local}
|u_i - v_i| \le |u_i| - |v_i| + 2|v_i|\sin\left(\frac{\theta}{2}\right)= |u_i| - c_\theta^{-1}|v_i|
\end{equation}

We now proceed via backward induction. For any time slice $i \in \{1, \dots, N\}$, our inductive hypothesis states
\begin{equation}\label{eq:induction_hyp}
\sum_{j=i}^{N-1} |v_j| \le c_{\theta} \sum_{j=i}^{N-1} |\gamma^*_{j+1} - \gamma^*_j| + c_{\theta} |\gamma_i - \gamma^*_i|
\end{equation}
If the induction hypothesis holds, evaluating at $i = 0$, and using $\gamma_0 = \gamma^*_0$, we obtain:
\[
\bar{A}(\gamma) = \sum_{j=0}^{N-1} |v_j| \le {c_{\theta}} \sum_{j=0}^{N-1} |\gamma^*_{j+1} - \gamma^*_j| = c_{\theta}  \bar{A}(\gamma^*),
\]
and the statement follows.

\textbf{Base Case:} For $i = N$, both sides of~\eqref{eq:induction_hyp} are $0$.

\textbf{Inductive Step:} Assume~\eqref{eq:induction_hyp} holds for time step $i+1$. \\
Since $|\gamma_{i+1} - \gamma^*_{i+1}| = |u_i - v_i|$, we apply~\eqref{eq:yao_local} and the definition of $c_\theta$ from \eqref{theta} to obtain:
\begin{align*}
|v_i| + c_{\theta}|\gamma_{i+1} - \gamma^*_{i+1}| \le & \ |v_i| + c_{\theta} \left( |u_i| - c_\theta^{-1}|v_i|  \right) \\
=& \ c_{\theta}|u_i| \\
=& \ c_{\theta}|(\gamma^*_{i+1} - \gamma^*_i) + (\gamma^*_i - \gamma_i)| \\
\leq&  \ c_{\theta}|\gamma^*_{i+1} - \gamma^*_i| + c_{\theta}|\gamma_i - \gamma^*_i|. 
\end{align*}
Then by the induction hypothesis, 
\begin{align*}
\sum_{j=i}^{N-1} |v_j| \le & \ c_{\theta}\sum_{j=i+1}^{N-1} |\gamma^*_{j+1} - \gamma^*_j| + |v_i| + c_{\theta} |\gamma_{i+1} - \gamma^*_{i+1}| \\
\leq & \ c_{\theta}\sum_{j=i+1}^{N-1} |\gamma^*_{j+1} - \gamma^*_j| + c_{\theta}|\gamma^*_{i+1} - \gamma^*_i| + c_{\theta}|\gamma_i - \gamma^*_i|
\\ = & \ c_{\theta}\sum_{j=i}^{N-1} |\gamma^*_{j+1} - \gamma^*_j| + c_{\theta}|\gamma_i - \gamma^*_i|.
\end{align*}

This completes the induction.

\end{proof}

\begin{lemma}\label{lem: covering}
For any $d \ge 1$, the unit sphere $\Sph^{d-1} \subset \mathbb{R}^d$ can be covered by at most~$7^d$ spherical caps of aperture $\pi/4$.
\end{lemma}

\begin{proof}
Let $\mathcal{N} \subset \Sph^{d-1}$ be a maximal set of points such that the angle between any two distinct points in $\mathcal{N}$ is strictly greater than $\pi/8$. 

By the maximality of $\mathcal{N}$, any point $x \in \Sph^{d-1}$ must have an angle of at most $\pi/4$ with at least one point $y \in \mathcal{N}$ (otherwise, $x$ could be added to $\mathcal{N}$, contradicting its maximality). Therefore, the spherical caps of angle $\pi/8$ (aperture $\pi/4$) centered at the points in $\mathcal{N}$ completely cover $\Sph^{d-1}$.

For any distinct $y_1, y_2 \in \mathcal{N}$, their distance is:
\[
|y_1 - y_2| = \sqrt{|y_1|^2 + |y_2|^2 - 2\langle y_1, y_2 \rangle} > \sqrt{2 - 2\cos(\pi/8)} \approx 0.390.
\]
Let $\eps =  \sqrt{2 - 2\cos(\pi/8)} $. Consider open Euclidean balls of radius $\eps/2$ centered at each point in $\mathcal{N}$. Due to the display above, these balls are pairwise disjoint. 
All these balls are contained within a single ball centered at the origin with radius $1 + \eps/2$. Comparing volumes yields:
\[
|\mathcal{N}| \cdot m\left(B\left(0, \frac{\eps}{2}\right)\right) \le m\left(B\left(0, 1 + \frac{\eps}{2}\right)\right).
\]
This implies
\[
|\mathcal{N}| \le \left( \frac{1 + \eps/2}{\eps/2} \right)^d = \left( 1 + \frac{2}{\eps} \right)^d\le 7^d,
\]
and the proof is complete.
\end{proof}

We recall that $\Gamma^n_\omega$ was defined in Section~\ref{sec:Intro} as the random collection of paths starting at $0$ and passing through Poissonian points.
\begin{lemma}\label{lem:abs}
Let $C,\eps>0$. Then there exists a $\delta > 0$ such that the following holds almost surely: for all sufficiently large $n$,
if
\begin{equation}\label{eq: bigbad}
\gamma\in\Gamma^n_{\omega}
\quad\text{and}\quad
\frac1n\sum_{k=0}^{n-1}
\mathbf 1_{\Delta_k\gamma \in S_{\delta}}
\geq
\eps,
\end{equation}
then
\begin{equation}\label{eq: biglength}
\bar A^n(\gamma)
\geq
Cn.
\end{equation}
\end{lemma}

\begin{proof}
Let's fix $n$, assume \eqref{eq: bigbad}, and set $R = \frac{2C}{\eps}$. If \[\frac1n\sum_{k=0}^{n-1}
\mathbf 1_{\Delta_k\gamma \in S_{\delta}} \mathbf{1}_{|\Delta_k\gamma| > R} \geq \frac{\eps}{2}, \] then~\eqref{eq: biglength} follows immediately, so we assume
\begin{equation}\label{eq: contradiction}
\frac1n\sum_{k=0}^{n-1}
\mathbf 1_{\Delta_k\gamma \in S_{\delta}} \mathbf{1}_{|\Delta_k\gamma| \leq R} \geq \frac{\eps}{2}.
\end{equation}
\newline
For each $k\ge0$, we define the random set
\begin{multline*}
D_{k, l}=
\Big\{
 (k,y_0)\in \pi\omega \cup \{(0,0)\} \ :  \exists(k+1,y_1), \dots, (k+l, y_l) \in \pi\omega, 
 \\
  y_{m+1} - y_m\in S_{\delta } \cap B(0, R)\  \textrm{for all}\ m=0,1,\ldots,l-1
\Big\}.
\end{multline*}
Notice that $D_{k, l+1} \subset D_{k, l}.$ 
For $I = (i_0, i_1, \dots, i_k), $ where $0 = i_0 < i_1 < i_2 < \dots < i_k \leq n$ and $\gamma \in \Gamma^n_{\omega}$, we define 
\begin{equation}
    \tilde{A}_I^n(\gamma) = \sum_{j=0}^{k-1}|\gamma_{i_{j+1}}-\gamma_{i_j}|.
\end{equation}
Note for all $I$ and $\gamma$, $\tilde{A}_I^n(\gamma) \leq \bar{A}^n(\gamma)$. Then for any path $\gamma$ satisfying condition~\eqref{eq: contradiction} but not 
satisfying~\eqref{eq: biglength}, 
there are $I = (i_0, i_1, \dots, i_k),$
$L = (l_0, l_1, \dots, l_k)$, $\beta>0$ such that 
\begin{equation}\label{eq: simple}
    \tilde{A}_I^n(\gamma) \leq Cn;
\end{equation}
\begin{equation}
\label{eq: constraint}
\begin{gathered}
    i_k + l_k \leq n -1, \\
    l_j \geq 1 \ \text{for $j\geq 1$}, \\ 
    i_{j+1} - (i_j + l_j) >0, \\ 
    \sum_{j=0}^k l_j = \lceil n\frac\eps2 \rceil -1 := \beta n, \\  
    \gamma_{i_j} \in D_{i_j, l_j},
\end{gathered}
\end{equation} 
and  
\[
\Delta_{i_j+m} \gamma \in S_{\delta},\quad |\Delta_{i_j +m}\gamma| \leq R,\quad \textrm{for all}\ m=0,1,\ldots,l_j-1.
\]
In fact, this provides a bijection between a choice of $\beta n $ indices from $0$ to $n-2$ and a choice of $I$ and $L$.

Now for fixed $I, L$ (which determines $k$), let $\gamma^*_{I, L}$ minimize $\tilde{A}^n_I$ among paths $\gamma \in \Gamma^n_{\omega}$ which have $\gamma_{i_j} \in D_{i_j, l_j}$ for all $0 \leq j \leq k$. By Lemma \ref{lem: covering}, we can cover $\Sph^{d-1}$ by $7^d$ spherical caps $\{K_m\}_{m=1}^{7^d}$ of aperture $ \frac{\pi}{4}$. By Lemma \ref{Yao}, if $\tilde{A}^n_I(\gamma^*_{I, L}) < Cn$, then there is some path $\gamma \in \Gamma^n_{\omega}$ such that for $0 \leq j \leq k-1$, $\gamma_{i_{j+1}}$ minimizes $|\gamma_{i_{j+1}}- \gamma_{i_j}|$ among points in $D_{i_{j+1}, l_{j+1}}$ in the direction of some cap, and 
\begin{equation}
    \tilde{A}^n_I(\gamma) < nc_{\theta}\cdot C.
\end{equation} 
There are $(7^d)^k$ possible paths for this $\gamma$, which can be indexed by $\sigma \in \{1, ..., 7^d\}^k, $ i.e., $\gamma^{I, L, \sigma} \in \Gamma^n_{\omega}$ is the path such that $\gamma^{I, L, \sigma}_{i_j}$ minimizes $|\gamma^{I, L, \sigma}_{i_{j+1}}- \gamma^{I, L, \sigma}_{i_j}|$ among points in $D_{i_{j+1}, l_{j+1}}$ in the direction of the cap $K_{\sigma_{j+1}}$.

Now we fix $I, L, \sigma$. For $0 \leq j \leq k-1$, the distribution of the random variable $|\gamma^{I, L, \sigma}_{i_{j+1}}-\gamma^{I, L, \sigma}_{i_j}|$ is equal to the distribution of the Poisson point $y$ closest to the origin among those contained in the cone formed by some fixed spherical cap~$K$ of aperture $\pi/4$ such that $(0, y) \in D_{0, l_{j+1}}$. In addition, the random variables $(|\gamma^{I, L, \sigma}_{i_{j+1}}-\gamma^{I, L, \sigma}_{i_j}|)_{j=0,\ldots,k-1}$ are mutually independent and independent of $\{(0, 0)\in D_{0, l_0}\}$, since the windows
$\{i_j,\ldots,i_j+l_j\}$ do not overlap and $\gamma_{i_j}$ is measurable with respect to the preceding layers.
Define 
\[
B_l^b := \Big\{(0, y) \in D_{0, l}: 0<|y| \leq b, y \in \R_+K\Big\}. 
\]
where $\R_+K=\{ry:\ r\ge0,\ y\in K\}$.
By Campbell's formula,
\begin{align*}
\mathbb E|B_l^b|=& \int_{B(0, b) \cap\, \R_+K}\mathbb P
\big((0,x)\in B_l^b\, \big|\, (0,x)\in\pi \omega\big)\,\dd x \\
=& \int_{B(0, b) \cap\, \R_+K}\mathbb P
\big((0,x)\in D_{0, l}\, \big|\, (0,x)\in\pi \omega\big)\,\dd x.
\end{align*}
By translation invariance,
\[
\Pp\big((0,x)\in B_l^b\, \big|\, (0,x)\in\pi \omega\big)=\Pp\big((0,y_0)\in B_l^b\, \big|\, (0,y_0)\in\pi \omega \big)
\]
for every $x\in\mathbb R^d$ and some fixed $y_0$.

Let $N$ denote the number of sequences $
(x_1,\ldots,x_{l}),$ $(i, x_i)\in\pi\omega,$
such that
\[
x_1 - y_0\in S_{\delta } \cap B(0, R), \qquad x_{i+1}-x_i\in S_{\delta } \cap B(0, R),\quad i=1,\ldots,l-1.
\]
Since $\{(0, y_0)\in D_{0, l}\}=\{N\ge1\}$,
\[
\Pp\big((0, y_0)\in D_{0, l}\, \big|\, (0,y_0)\in\pi \omega \big) \le \E [N \,|\,(0,y_0)\in\pi \omega].
\]
Due to the multivariate Campbell formula and the independence property of the Poisson process, the  right-hand side equals
\[
\int_{(\mathbb R^d)^l} \mathbf 1_{\{x_1- y_0\, \in\, S_{\delta } \cap B(0, R)\}} \prod_{i=1}^{l-1} \mathbf 1_{\{x_{i+1}-x_i\,\in\, S_{\delta } \cap B(0, R)\}}\,\dd x_1\cdots \dd x_{l}=
m(S_{\delta } \cap B(0, R))^l.
\]
Therefore,
\[
\E |B_l^b|\le m(B(0, b) \cap\, \R_+K)\,(m(S_{\delta } \cap B(0, R)))^l \leq
V_db^d(m(S_{\delta } \cap B(0, R)))^l.
\]
Then we obtain
\begin{equation*}
    \mathbb{P}\big\{|\gamma^{I,L,\sigma}_{i_{j+1}}-\gamma^{I,L,\sigma}_{i_j}| < b\big\}
    = \mathbb{P}\big\{|B_{l_{j+1}}^b| \geq 1\big\} \leq V_db^d(m(S_{\delta } \cap B(0, R)))^{l_{j+1}}.
\end{equation*}
A similar calculation gives 
\[
\mathbb{P}\{(0, 0)\in D_{0, l_0}\} \leq (m(S_{\delta } \cap B(0, R)))^{l_{0}}.
\]
Then, by Strassen's theorem \cite{Strassen:MR0177430}, after possibly enlarging the underlying probability space, there exists mutually independent  random variables \(U_j\) independent of $\{(0, 0)\in D_{0, l_0}\}$ such that
\[
|\gamma^{I,L,\sigma}_{i_{j+1}}-\gamma^{I,L,\sigma}_{i_j}| \ge U_j
\]
almost surely, where \(U_j\) has density
\[
f_{U_j}(b)
= V_d d b^{d-1}
m\bigl(S_{\delta}\cap B(0,R)\bigr)^{l_{j+1}}
\mathbf{1}_{\left[0,\,(V_d m(S_{\delta}\cap B(0,R))^{l_{j+1}})^{-1/d}\right]}(b).
\]
By applying Dirichlet's integral formula over the simplex $\left\{ \sum_{j=0}^{k-1} b_j < n c_{\theta} C \right\}$ weighted by the joint density of the independent random variables $U_j$, we obtain: 
\begin{equation*}
\begin{split}
\mathbb{P}\Bigg\{ \sum_{j=0}^{k-1} U_j < n c_{\theta} C \Bigg\} \leq
\int_{\sum_{j=0}^{k-1} b_j < n c_{\theta} C} \prod_{j=0}^{k-1}V_ddb_j^{d-1}m\big(S_{\delta} \cap B(0, R)\big)^{l_{j+1}}\mathrm{d}b_0...\mathrm{d}b_{k-1} \\
= \int_{\sum_{j=0}^{k-1} b_j < 1} \prod_{j=0}^{k-1}V_dd(b_jnc_{\theta}C)^{d-1}m\big(S_{\delta} \cap B(0, R)\big)^{l_{j+1}}(nc_{\theta}C\mathrm{d}b_0)...(nc_{\theta}C\mathrm{d}b_{k-1})
\\ \leq \frac{\Gamma(d)^kd^k}{\Gamma(kd+1)} (n c_{\theta} C)^{kd} \prod_{j=0}^{k-1} V_d m\big(S_{\delta} \cap B(0, R)\big)^{l_{j+1}}.
\end{split}
\end{equation*}
Using \eqref{eq: sausage} and  the constraint $\sum_{j=0}^{k} l_{j} = n\beta $ 
to estimate the product, we get
\begin{equation}
\label{eq: single_bound}
\mathbb{P}\big\{(0, 0)\in D_{0, l_0}\big\}\mathbb{P}\Bigg\{ \sum_{j=0}^{k-1} U_j < n c_{\theta} C \Bigg\}\le
\frac{\Gamma(d)^kd^k}{\Gamma(kd+1)} (n c_{\theta} C)^{kd} V_d^k (k_R\delta^{\eta_R})^{n\beta}.
\end{equation}

Now, let $P_{n}$ denote the total probability that there exists at least one path $\gamma \in \Gamma^n_{\omega}$ satisfying \eqref{eq: contradiction} but not \eqref{eq: biglength}. Our plan is to check that $P_n$ decays exponentially. Then the Borel--Cantelli lemma will imply that such paths do not exist for sufficiently large $n$, and the claim of our lemma will follow.

To that end, we apply a union bound over all valid choices of index sets~$I$, continuous block lengths $L$, and directional profiles $\sigma$. 
The total number of configurations for $(I, L)$ is bounded by $\binom{n}{\beta n}$, and for a fixed $k \leq \beta n$, the number of directional cap assignments $\sigma$ is $(7^d)^k \leq 7^{d\beta n}$. Combining these with \eqref{eq: single_bound} yields:
\begin{align*}
P_{n} &\leq \sum_{I, L, \sigma}\mathbb{P}\big\{(0, 0)\in D_{0, l_0}\big\} \mathbb{P}\Bigg\{ \sum_{j=0}^{k-1} U_j < n c_{\theta} C \Bigg\} 
\leq J_1(n) J_2(n) J_3(n) (k_R\delta^{\eta_R})^{n\beta}, 
\end{align*}
where 
\begin{align*}
    J_1(n) =  \binom{n}{n\beta}, \quad
    J_2(n) =  7^{dn\beta}, \quad
    J_3(n) := \max_{1 \leq k \leq n\beta} \left[ \frac{\big( (n c_{\theta} C  )^{d} \Gamma(d) dV_d) \big)^k}{\Gamma(kd+1)} \right].
\end{align*}
To ensure that $P_n$ decays exponentially, it suffices to 
show that $J_i(n)$, $i=1,2,3,$ grow at most exponentially and choose $\delta$ sufficiently small.

For $J_1(n)$, we use the bound 
\begin{equation*}
    \binom{n}{k} \leq \frac{n^k}{k!} \leq \frac{n^k}{(k/e)^k} \leq \left(\frac{en}{k}\right)^k,
\end{equation*}
where $(k/e)^k \leq k!$ follows from the evaluation of the Taylor series for $e^x$ at $x=k$. Thus, we obtain the desired exponential bound
\begin{equation*}
J_1(n) = \binom{n}{n\beta} \leq \left(\frac{en}{n\beta}\right)^{n\beta} = \left(\frac{e}{\beta}\right)^{\beta n}.
\end{equation*}

For $J_2(n)$, the exponential growth is obvious.

For $J_3(n)$, let $A = (c_{\theta} C)^d \Gamma(d) d V_d$. We need to bound the maximum over $k$ of $\frac{(A n^d)^k}{\Gamma(kd+1)}$. Using the bound $\Gamma(kd+1) = (kd)! \geq (kd/e)^{kd}$, we get:
\begin{equation*}
\frac{(A n^d)^k}{\Gamma(kd+1)} \leq \frac{A^k n^{dk}}{(kd/e)^{kd}} = \left( \frac{A e^d}{d^d} \left(\frac{n}{k}\right)^d \right)^k\le \left( \frac{A e^d}{d^d}\vee 1\right)^n \left(\frac{n}{k}\right)^{dk}
\end{equation*}
Thus, it suffices to show that the function 
\[
f_n(x)=\ln \left(\frac{n}{x}\right)^{dx}=d( x\ln n- x\ln x)  ,\quad 0<x\le n,
\]
admits a linear uniform bound in $n$. A simple calculation implies that its maximum is
attained at $ne^{-1}$ and equals $dne^{-1}$. This establishes the desired linear bound on $f_n$ and the exponential growth of $J_3(n)$ completing the proof of the lemma.
\end{proof}

We recall that for $v\in\R^d$, the path $\gamma_B^n(v)$ is defined via~\Cref{prop:min}. It provides the minimum in the definition of $B^n(v)$ given in \eqref{eqn:Bnx}.

\begin{corollary}\label{cor: key}
Almost surely with respect to $\Pp$,
\[
\lim_{\delta \downarrow 0}\limsup_{n \to \infty} \frac{1}{n}\sum_{k=0}^{n-1}\mathbf{1}_{\Delta_k\gamma_B^n(v) + v\, \in S_{\delta} } = 0.
\]

\end{corollary}
To derive this corollary, we will use the following lemma.
\begin{lemma}\label{lem:Lambda growth}
For all $v\in\R^d$, we have
$
\Lambda(v)\ge a|v|+b. 
$
\end{lemma}

\begin{proof}
Our assumption~\eqref{eq: growth of L} 
implies that  for all $\omega$ and for all paths $\gamma\in\Gamma^n(nv)$, we have
\[
A^n_\omega(\gamma)\ge \sum_{k=0}^{n-1} (a|\Delta_k\gamma|+b)\ge n(a|v| +b), 
\]
and our claim follows.
\end{proof}

\begin{proof}[Proof of Corollary~\ref{cor: key}]
    Let $\eps>0$. Lemma~\ref{lem:Lambda growth} implies that 
    $C = \frac{\Lambda(v)+1-b}{a}$ is positive. 
    Using \eqref{eq: growth of L} and applying Lemma \ref{lem:abs} 
    with this choice of~$C$,
    we  
    obtain that there is $\delta=\delta(\eps)>0$ such that with probability 1, for sufficiently large $n$, for all $\gamma \in \Gamma^n_{\omega}$, $\frac{1}{n}\sum_{k=0}^{n-1}\mathbf{1}_{\Delta_k\gamma \in S_{\delta}} > \eps$ implies 
    \[
    \frac{1}{n}\sum_{k=0}^{n-1}L(\Delta_k\gamma)  \ge \frac{1}{n}(a\bar A^n(\gamma)+bn)
    \ge aC+b
    \ge \Lambda(v) + 1.
    \]
    But we know that with probability 1, $\frac{1}{n}\sum_{k=0}^{n-1}L(\Delta_k\gamma_A^n(nv)) < \Lambda(v) + 1$ for sufficiently large $n$, where $\gamma_A^n(nv)$ is the minimizer for the action in the untransformed environment described in \Cref{prop:min}. Thus,
    \[
\limsup_{n \to \infty} \frac{1}{n}\sum_{k=0}^{n-1}\mathbf{1}_{\Delta_k\gamma_A^n(nv) \in S_{\delta(\eps)} } \le \eps.
\]
 Rewriting this in terms of $\Xi_v^*\omega$, letting $\eps\downarrow 0$, and noting $S_{\delta_0} \subseteq S_{\delta_1}$ whenever $\delta_0 \leq \delta_1$, we  complete the proof.
\end{proof}

\begin{lemma}\label{lem: K}
Let
\[
K^n(\delta)=\frac{1}{n}\sum_{k=0}^{n-1}K(\Delta_k\gamma_B^n(v) + v)\mathbf{1}_{\Delta_k\gamma_B^n(v) + v\, \in S_{\delta} }.
\]
Then, almost surely with respect to $\Pp$,
\[
\lim_{\delta  \downarrow 0}\limsup_{n \to \infty} K^n(\delta) = 0.
\]

\end{lemma}

\begin{proof}
Set $f^n_k(v)=\Delta_k\gamma^n_B(v) + v.$ By Hölder's inequality with exponents $1/(1-\theta)$ and $1/\theta$ and \eqref{eq: K} (it suffices to consider that condition for $\theta\ne 0$),
\begin{multline*}
    K^n(\delta)=\frac{1}{n}\sum_{k=0}^{n-1}K(f^n_k(v))\mathbf{1}_{f^n_k(v) \in S_{\Dwv} } 
    \leq  \frac{1}{n}\sum_{k=0}^{n-1}\big(c_0 + c_1|L(f^n_k(v)) |^{\theta}\big)\mathbf{1}_{f^n_k(v) \in S_{\Dwv} } 
    \\
    \leq 
     \frac{c_0}{n}\sum_{k=0}^{n-1}\mathbf{1}_{f^n_k(v) \in S_{\Dwv} } 
     +  
    c_1\Big(\frac{1}{n}\sum_{k=0}^{n-1}\mathbf{1}_{f^n_k(v)\in S_{\Dwv}}\Big)^{1-\theta}\Big(\frac{1}{n}\sum_{k=0}^{n-1}|L(f^n_k(v))|\mathbf{1}_{f^n_k(v)\in S_{\Dwv}}\Big)^{\theta}.
\end{multline*}
We know with probability 1, 
\begin{align*}
\frac{1}{n}\sum_{k=0}^{n-1}|L(f^n_k(v))| \mathbf{1}_{f^n_k(v) \in S_{\Dwv}} \leq \frac{1}{n}\sum_{k=0}^{n-1}L(f^n_k(v)) + 2|b| \leq \Lambda(v) + 1 + 2|b|
\end{align*}
for sufficiently large $n$. The result follows from Corollary \ref{cor: key}. 
\end{proof}

\begin{lemma}\label{lem: H}
Let 
\[
H^n(\delta)=\frac{1}{n}\sum_{k=0}^{n-1}\delta ^\alpha H(\Delta_k\gamma_B^n(v) + v)\ONE _{\Delta_k\gamma^n_B(v)+v\, \notin S_{\delta}}.
\]
Then,
almost surely with respect to $\Pp$,
\begin{equation}
\label{eq:I}
 \lim_{\delta \downarrow 0}\limsup_{n \to \infty} H^n(\delta) = 0.
\end{equation}

\end{lemma}
\begin{proof}
    Set $f^n_k(v)=\Delta_k\gamma^n_B(v) + v.$
    Then, recalling how $r_0>0$ is used in our assumptions~\eqref{eq: H} and~\eqref{eq: H local}, for any $\eps\in(0,r_0)$ we can write 
\begin{equation*}
H^n(\Dwv)
= H^n_1(\Dwv)+H^n_2(\Dwv,\eps)+H^n_3(\Dwv,\eps),
\end{equation*}
where
\begin{align*}
    H^n_1(\Dwv) &= \frac{1}{n}\sum_{k=0}^{n-1}\Dwv^\alpha H(f^n_k(v))
\ONE_{f^n_k(v)\notin S_{\Dwv}}
\ONE_{f^n_k(v)\notin S_{r_0}},
\\ 
    H^n_2(\Dwv,\eps) &= \frac{1}{n}\sum_{k=0}^{n-1}\Dwv^\alpha H(f^n_k(v))
\ONE_{f^n_k(v)\notin S_{\Dwv}}
\ONE_{f^n_k(v)\in S_{r_0}}\ONE_{f^n_k(v)\in S_{\eps}},
\\
    H^n_3(\Dwv,\eps) &= \frac{1}{n}\sum_{k=0}^{n-1}\Dwv^\alpha H(f^n_k(v))
\ONE_{f^n_k(v)\notin S_{\Dwv}}
\ONE_{f^n_k(v)\in S_{r_0}}\ONE_{f^n_k(v)\notin S_{\eps}}.
\end{align*}
We can use~\eqref{eq: H} to obtain for sufficiently large $n$,
\begin{align*}
    H^n_1(\Dwv) &\leq \frac{\Dwv^\alpha}{n}\sum_{k=0}^{n-1}
\big(c_0+c_1|L(f_k^n(v))|\big)
\\ & \leq \frac{\Dwv^\alpha}{n}\sum_{k=0}^{n-1}
\big(c_0+c_1(L(f_k^n(v))+2|b|)\big)
\\
&\leq  \Dwv^\alpha(c_0 + c_1(\Lambda(v) + 1+2|b|)),
\end{align*}
so
\[
 \lim_{\Dwv\downarrow 0} \limsup_{n \to \infty}H^n_1(\Dwv)  = 0.
\]
Also, for fixed $\eps$, we use \eqref{eq: H local} to obtain
\begin{align*}
    H^n_3(\Dwv,\eps) \leq \frac{\Dwv^\alpha}{\eps^\alpha}\cdot\frac{1}{n}\sum_{k=0}^{n-1}\eps^\alpha H(f^n_k(v))
\ONE_{f^n_k(v)\in S_{r_0}}\ONE_{f^n_k(v)\notin S_{\eps}}\le \frac{c_2\Dwv^\alpha}{\eps^\alpha},
\end{align*}
so that
\[
 \lim_{\Dwv\downarrow 0}\limsup_{n \to \infty}H^n_3 (\Dwv,\eps) = 0.
\]
Thus, for any $\eps>0$,
\begin{equation}
\label{eq:estim_In}
\limsup_{\delta  \downarrow 0} \limsup_{n\to\infty} H^n(\delta)\le \limsup_{\Dwv\downarrow 0} \limsup_{n\to\infty} H^n_2(\Dwv,\eps).
\end{equation}
Due to~\eqref{eq: H local}, 
\begin{equation*}
H^n_2(\Dwv,\eps) \leq \frac{c_2}{n}\sum_{k=0}^{n-1}\ONE_{f^n_k(v)\in S_{\eps}}. 
\end{equation*}
Since the right-hand side of this estimate does not depend on $\Dwv$, for any fixed $\eps>0$ we can 
continue~\eqref{eq:estim_In} as
\begin{equation*}
\limsup_{\delta \downarrow 0} \limsup_{n\to\infty} H^n(\delta)\le \limsup_{n\to\infty} \frac{c_2}{n}\sum_{k=0}^{n-1}\ONE_{f^n_k(v)\in S_{\eps}}.
\end{equation*}
Taking $\eps\downarrow 0$ in this estimate and using Corollary~\ref{cor: key}, we complete the proof.
\end{proof}

\medskip

\begin{proof}[Proof of theorem] 
For a convex function, the $C^1$ property is equivalent to differentiability at every point  (see \cite[Theorem 25.5]{Rockafellar1970}). So let's fix an arbitrary $v\in\R^d$ and prove differentiability of $\Lambda$ at $v$ with the help of Lemma~\ref{lem: conv}
on approximate linear domination.

For $i\in\{1, \dots, d +1\}$, 
it follows from \eqref{eq: linear} that for any path $\gamma \in \Gamma^n(0)$,
\begin{align*}
\frac{1}{n} B^n_\omega(v+\eps u_i, \gamma) &\le \frac{1}{n} B^n_\omega(v, \gamma) + \frac{1}{n}\sum_{k=0}^{n-1} \langle M(\Delta_k\gamma + v),\eps u_i \rangle\\
&\quad + \frac{\eps^{1+\alpha}}{n} \sum_{k=0}^{n-1}  H(\Delta_k\gamma + v)\mathbf{1}_{\Delta_k\gamma+ v\,\notin S_{\eps}} \\ 
&\quad + \frac{\eps}{n} \sum_{k=0}^{n-1} K(\Delta_k\gamma + v )\mathbf{1}_{\Delta_k\gamma+ v\, \in S_{\eps}}.
\end{align*}

Let $\gamma^n_{A,\Xi_v^*\omega}(nv)$ be a measurable sequence of minimizers the existence of which follows from Proposition~\ref{prop:min} and the invariance of $\Pp$ under $\Xi_v^*$.
Then 
$\gamma^n_{B,\omega}(v)=\Xi_v^{-1}\gamma^n_{A,\Xi^*_v\omega}(nv)$ is a measurable sequence of minimizers of the action~$B^n$.  

For brevity, we denote these paths by $\gamma^n_B(v)$. Plugging them
into the above display, using the bound $B^n_{\omega}(v+\eps u_i) \le B^n_{\omega}(v+\eps u_i, \gamma^n_B(v))$ and the identity $B^n_{\omega}(v) = B^n_{\omega}(v, \gamma^n_B(v))$, we obtain
\begin{align}\notag
\frac{1}{n} B^n_{\omega}(v+\eps u_i, \gamma^n_B(v)) &\le \frac{1}{n} B^n_{\omega}(v, \gamma_B^n(v)) + \frac{1}{n}\sum_{k=0}^{n-1} \langle M(\Delta_k\gamma_B^n(v) + v), \eps u_i \rangle \\
\label{eq:checking-linear-domination}
&\quad + \frac{\eps^{1+\alpha}}{n} \sum_{k=0}^{n-1}  H(\Delta_k\gamma_B^n(v) + v)\mathbf{1}_{\Delta_k\gamma_B^n(v)+ v\,\notin S_{\eps}} \\ \notag
&\quad + \frac{\eps}{n} \sum_{k=0}^{n-1} K(\Delta_k\gamma_B^n(v) + v )\mathbf{1}_{\Delta_k\gamma_B^n(v)+ v\, \in  S_{\eps}}.
\end{align}
Now we can apply Lemma~\ref{lem: conv} 
\begin{align*}
    x_0 &= v, \\
    x_{m}^{(i)} &= v + \frac{u_i}{m}, \\
    f_n(x) &= \frac{1}{n}B^n(x), \\
    f &= \Lambda, \\
    \xi_n &= \frac{1}{n}\sum_{k=0}^{n-1}M(\Delta_k\gamma_B^n(v) + v),
\end{align*}
and 
\begin{multline*}
    h_{n,m}^{(i)} = \frac{(\frac{1}{m})^{1+\alpha}}{n} \sum_{k=0}^{n-1}  H(\Delta_k\gamma_B^n(v) + v)\mathbf{1}_{\Delta_k\gamma_B^n(v)+ v\,\notin S_{1/m}}  \\ + \frac{1}{mn} \sum_{k=0}^{n-1} K(\Delta_k\gamma_B^n(v) + v )\mathbf{1}_{\Delta_k\gamma_B^n(v)+ v\,\in S_{1/m}}.
\end{multline*}

The first three conditions of Lemma~\ref{lem: conv} are straightforward consequences of these definitions and the assumptions of the theorem. Conditions~\eqref{eq:f_n-conv-to-f-at-x_m} 
and~\eqref{eq:f_n-conv-to-f-at-x_0} of Lemma \ref{lem: conv} follow from~Theorem~\ref{thm: shape} and Remark~\ref{rem:Lambda-via-B}. 
Condition~\eqref{limit} of Lemma \ref{lem: conv} holds almost surely due to Lemmas \ref{lem: K} and \ref{lem: H}. Condition \eqref{eq:lemma33-assumption} follows from~\eqref{eq:checking-linear-domination} with $\eps=1/m$. Thus Lemma~\ref{lem: conv} implies
differentiability of $\Lambda$ at $v$. The almost sure relation~\eqref{eq:gradient} 
with $\omega$ replaced by $\Xi^*_v\omega$
follows 
from~\eqref{eq:lemma33-conclusion}. 
Since $\Pp$ is invariant under $\Xi^*_v$, \eqref{eq:gradient} also holds for almost all $\omega$.
Note that although uniqueness of minimizers may fail, our argument works for any measurable selection of the minimizers.
The proof is completed. 
\end{proof}

\section{Proof of Theorems on Limit Shapes}
\label{sec:limit-shape-proofs}
\subsection{Proof of Theorem~\ref{thm:footprint}}
\label{sec:footprint}
First, we note that $H_\Lambda=\partial E_\Lambda\cap H
= E_\Lambda\cap H$
is convex since $E_\Lambda$ and $H$ are convex.

Lemma~\ref{lem:Lambda growth} and assumption~\eqref{eq:L-nonnegative} imply that for all $u\in\Sph^{d-1}$ and all $t\in(0,\Lambda(0)^{-1})$, we can find $r_u(t)$ such that
\begin{equation}
\label{eq:defining r t}
t\Lambda\Big(\frac{r_u(t)u}{t}\Big)=1.    
\end{equation}
Then $\{(t,r_u(t)u): t\in(0,\Lambda(0)^{-1})\}\subset \partial E_\Lambda\setminus H$ and, defining 
$r_u(0)=\lim_{t\downarrow 0}r_u(t)$, we obtain 
\[
\{r\ge 0: (0,ru)\in H_\Lambda\}=[0,r_u(0)].
\]
Thus the theorem will follow once we prove
\begin{equation}
\label{eq: goal}
r_u(0)=\lambda(u)^{-1}\in[0,\infty),    
\end{equation}
where $\lambda(u)$ was defined in \eqref{eq:linear rate lambda u}.
Rewriting~\eqref{eq:linear rate lambda u} as
\[
\Lambda(su)=\frac{s}{\lambda(u)^{-1}+o(1)},\quad s\to +\infty,
\]
and noting that \eqref{eq:defining r t} implies $\lim_{t\downarrow 0}(r_u(t)/t)=+\infty$,
we obtain from \eqref{eq:defining r t}
\[
t\frac{r_u(t)/t}{\lambda(u)^{-1}+o(1)}=1,
\]
implying \eqref{eq: goal} and completing the proof.
\subsection{Proof of Lemma~\ref{lem:Lambda-positive}}

Applying  Corollary~\ref{cor: key} of Section~\ref{sec:main_proofs} to paths and environments transformed under shear, we obtain
that there are $\delta,c>0$ such that if $n$ is sufficiently large, the action minimizer $\gamma^n_A(nv)$ satisfies  
\[
\sum_{k=0}^{n-1}\ONE_{\gamma^n_{A}(nv)_{k+1}-\gamma^n_A(nv)_{k}\, \notin\, S_\delta}>cn,
\]
so that $\Lambda(v)\ge cL_*(\delta)>0$, and the proof is complete.

\subsection{Proof of Theorem~\ref{thm: d>1}} 
\label{sec:d>1}
Throughout this section, $d\ge 2$, $L(v)=|v|$, $\nu=\delta_0$.
\begin{lemma}\label{lem: lowerBound}
    For all $v\ne0$, $\Lambda(v) > |v|. $
\end{lemma}
\begin{proof}
    By the triangle inequality, we know that $\Lambda(v) \geq |v|.$ Assume for the sake of contradiction that there exists some $v_0$ with $\Lambda(v_0) = |v_0|$. Since $\Lambda$ is radially symmetric, we may assume $v_0 = (x_0, 0, ..., 0), $ where $x_0 = |v_0|$. For each $n$, let $\gamma^n = (\gamma^n_0 =0, \gamma^n_1, ..., \gamma^n_n = nv_0)$ be a minimizing path to \eqref{eqn:Anx}, where $\gamma_k^n = (\gamma_k^{n, 1}, ..., \gamma_k^{n, d})$. By our assumption, $\sum_{k=0}^{n-1} |\gamma^n_{k+1} - \gamma^n_k| = x_0n + o(n)$ a.s. 

    Now define the map $T$ acting on paths $\gamma $ by 
    \[
    (T\gamma)_k = \Big(\frac{1}{2}\gamma^1_k, 2\gamma^2_k, \gamma^3_k, ..., \gamma^d_k\Big),
    \]
    and define the map $T^*$ acting on the environment, $T^*: \Omega \to \Omega$ by 
    \[
    T^*\omega = \Big\{\Big(k, \Big(\frac{x_1}{2}, 2x_2, x_3, ..., x_d\Big), \varphi\Big) : (k, (x_1, ..., x_d), \varphi) \in \omega\Big\}.
    \]
    By our assumptions on $\omega$, we have 
    $T^*\omega \overset{d}{=} \omega,$ 
    so
    \[
    A^n_{T^*\omega}(x) \overset{d}{=} A_\omega^n(x).
    \]
    The path $\overline\gamma^n = T\gamma^n$ belongs to $\Gamma^n_{T^*\omega}$ with endpoint $\overline\gamma^n_n = (n\frac{x_0}{2}, 0, ..., 0)$. We have 
    \begin{multline*}
        A^n_{T^*\omega}\Big(\Big(n\frac{x_0}{2}, 0, ..., 0\Big)\Big)  \leq  \sum_{k= 0}^{n-1}|\Delta_k\overline{\gamma}^n| \\
         = \sum_{k= 0}^{n-1}\Big( \frac{1}{4}(\Delta_k\gamma^{n, 1})^2 + 4(\Delta_k\gamma^{n, 2})^2 + (\Delta_k\gamma^{n, 3})^2 + ... + (\Delta_k\gamma^{n, d})^2 \Big)^{1/2} \\
        \leq \frac{1}{2}A^n_\omega((nx_0,0, ..., 0)) + 2\sum_{k=0}^{n-1} \Big((\Delta_k\gamma^{n, 2})^2 + ... + (\Delta_k\gamma^{n, d})^2 \Big)^{1/2}.
       \end{multline*}
       Then with probability 1,
       \begin{align*}
        \,& 2\sum_{k=0}^{n-1} \Big((\Delta_k\gamma^{n, 2})^2 + ... + (\Delta_k\gamma^{n, d})^2 \Big)^{1/2}  
        \\ 
        =\,&  2\sum_{k=0}^{n-1} \Big(|\Delta_k{\gamma}^n|^2 - |\Delta_k\gamma^{n, 1}|^2\Big)^{1/2}  \\
        = \, &   2\sum_{k=0}^{n-1} \Big((|\Delta_k{\gamma}^n| - |\Delta_k\gamma^{n, 1}|)(|\Delta_k{\gamma}^n| + |\Delta_k\gamma^{n, 1}|)\Big)^{1/2} \\ 
         \leq \, &  2 \bigg(\sum_{k=0}^{n-1} \big(|\Delta_k{\gamma}^n| - |\Delta_k\gamma^{n, 1}|\big)\sum_{k=0}^{n-1} \big(|\Delta_k{\gamma}^n| + |\Delta_k\gamma^{n, 1}| \big) \bigg)^{1/2} \\ 
        \leq \, & 2 \Big(A_\omega^n((nx_0, ..., 0)) - \sum_{k= 0}^{n-1}\Delta_k \gamma^{n, 1}\Big)^{1/2} (2A_\omega^n((nx_0, ..., 0)))^{1/2} \\
        = \, & 2^{3/2} \big(nx_0 + o(n) - x_0n\big)^{1/2} \big(nx_0 + o(n) \big)^{1/2} = o(n).
    \end{align*}
  Combining the last two displays, we obtain an almost sure estimate
    \[
    A^n_{T^*\omega}\Big(\Big(n\frac{x_0}{2}, 0, ..., 0\Big)\Big) \leq \frac{1}{2}A^n_\omega((nx_0,0, ..., 0)) + o(n) = n\frac{x_0}{2} + o(n)
    \]
  implying $\Lambda((\frac{x_0}{2}, 0, ..., 0)) = \frac{x_0}{2}$. Iteratively applying this argument and using the symmetry of $\Lambda$, we obtain
    \[
    \Lambda\Big(\Big(\pm\frac{x_0}{2^m}, 0, ..., 0\Big)\Big)=\frac{x_0}{2^m}, \quad m\in\N,
    \]
    but this implies that  the directional derivatives of $\Lambda$ along $(\pm1,0,\ldots,0)$ do not agree,
    contradicting the differentiability of $\Lambda$ at $0$. Hence $\Lambda(v) > |v|$. 
\end{proof}

    \begin{lemma}\label{lemma: lineargrowth}
    We have 
    \[
    \lim_{|v| \to \infty}(\Lambda(v) -|v|) = 0.
    \]

    \end{lemma}
    \begin{proof}
        Since $\Lambda(v)$ is rotationally symmetric, it suffices to consider $v = (x, 0, ..., 0)$. Consider the path $\gamma^n$ with $\gamma^n_0 = 0, \gamma^n_n = (nx, 0, ..., 0)$, and for $0 \leq k \leq n-2, $ $\gamma^n_k$ is defined recursively as follows: for $\gamma^n_k = (\gamma^{n, 1}_k, ..., \gamma^{n, d}_k)$, $\gamma^n_{k+1} = (\gamma^{n, 1}_{k+1}, ..., \gamma^{n, d}_{k+1})$ is chosen so that $(k+1, \gamma^n_{k+1})$ 
        is an element of $\pi\omega$ with the smallest value of $\gamma_{k+1}^{n, 1}$ satisfying the conditions
        \begin{itemize}
            \item $\gamma_{k+1}^{n, 1} \geq \gamma_{k}^{n, 1}$
            \item $\big((\Delta_k\gamma^{n, 2})^2 + ...+(\Delta_k\gamma^{n, d})^2\big)^{1/2} \leq R$. 
        \end{itemize}
        Then for $0 \leq k \leq n-2,$ $\Delta_k\gamma^n$ are i.i.d.\ random variables with first coordinate distributed exponentially with  parameter 
        $\lambda = V_{d-1} R^{d-1}$.
        Let's choose $R=R(x)$ so that $\lambda = \frac{1}{x}$, i.e. $\mathbb{E}[\Delta_k\gamma^n] = (x, 0, ..., 0)$ (all the components except the first one vanish by symmetry), and write 
        \begin{align*}
            \frac{1}{n}A^n_\omega((nx, 0, ..., 0))
        \le \frac{1}{n}\sum_{k=0}^{n-2}|\Delta_k\gamma^n| + \frac{1}{n}\Big|xn-\sum_{k=0}^{n-2}\Delta_k\gamma^n\Big|.
        \end{align*}
        By the law of large numbers, as $n \to \infty$, the first term converges almost surely to $\mathbb{E}[|\Delta_0\gamma^n|]$, and  the second term converges almost surely to $0$. Since 
         $|\Delta_0\gamma^n|\le \Delta_0\gamma^{n, 1} + R$,
        we have
        \begin{align*}
            \mathbb{E}[|\Delta_0\gamma^n|] \leq \mathbb{E}[\Delta_0\gamma^{n, 1}] + R = x +R.
        \end{align*}
        Therefore $\Lambda((x, 0, ..., 0)) \leq x+R$. As $x \to \infty$, $R \to 0$. Then \\ $\lim_{x \to \infty}\Lambda((x, 0, ..., 0))- x \leq 0$. Using radial symmetry of $\Lambda$ and that $\Lambda(v) \geq |v|$ implies the desired result.
    \end{proof}

    \begin{proof}[Proof of Theorem~\ref{thm: d>1}]
        By Lemma \ref{lemma: lineargrowth}, $\lambda(u) = 1$, so Theorem \ref{thm:footprint} implies 
        \[
        H_\Lambda=\{0\} \times \Big\{ru:\ 0\le r \le 1\Big\}.
        \]
        So $H_\Lambda$ is the closed $d$-dimensional unit ball. 
        The closure of the effective front contains the points $(t, x)$ satisfying $t\Lambda(x/t) =1$, along with the points $(0, x)$ satisfying $\lim_{t\downarrow 0}t\Lambda(x/t) = 1$. By Lemma \ref{lemma: lineargrowth}, the latter are the points $(0, x)$ with $|x| = 1$. Since the limit shape is radially symmetric, the closure of the effective front is exactly the set of points swept out by spatial rotations of the generating curve. 
        
        Suppose for the sake of contradiction that the generating curve was linear at $\big(0, e_1\big)$. Then there exists some $\kappa$ such that for  $t > 0$ sufficiently small, the points $(0, e_1) + t(1, \kappa e_1)$ are contained in the effective front. This implies $t\Lambda(e_1/t + \kappa e_1) = 1.$ Setting $v = e_1/t  + \kappa e_1$ gives $\Lambda(v) = |v| - \kappa$ for large enough $v$. By Lemma \ref{lemma: lineargrowth}, $\kappa $ must be $0,$ but this contradicts Lemma \ref{lem: lowerBound}. 
    \end{proof}

\section{Proof of Theorem~\ref{thm: shape}}
\label{sec:conv-to-shape}
The standard way to prove existence of shape functions is to invoke the subadditive ergodic theorem.
Subadditivity for the minimal action~$A$ writes as 
\begin{equation}
\label{eq:subadditivity}
A^{n+m}_\omega (x+y)\le A^{n}_\omega (x)+A^m_{\Theta^{n,x}\omega}(y),    
\end{equation}
where $\Theta^{n,x}:\Omega\to\Omega$ is a $\Pp$-preserving and ergodic shift mapping each Poisson point $(k,z,\varphi)$ to $(k-n,z-x,\varphi)$. 

This property is usually proved by noting that the action of the concatenation of two admissible paths sharing one endpoint is an upper bound for the optimal action between their remaining endpoints. In our case, the concatenation of two admissible paths is typically not admissible since it may contain a non-Poissonian point. However, subadditivity still holds for an auxiliary action $\tilde A$ approximating $A$. 

To define $\tilde A$,
we extend the Poisson point process to be defined on
$(\{0\}\cup\N)\times\R^d\times\R_+$.
Then, for any $\omega\in\Omega$ and any path $\gamma\in \Gamma^n_\omega$  we define 
$\tilde\gamma$ to be a path $\tilde\gamma_k:\{0,\ldots,n\}\to\R^d$ such that for all $k=1,\ldots,n-1$,
$\tilde\gamma_k=\gamma_k$, and $\tilde\gamma_0$ and $\tilde\gamma_n$ are the spatial components of Poissonian points in layers $0$ and $n$, closest to points $\gamma_0=0$ and $\gamma_n$, respectively. If there is more than one such a point, we choose the smallest with respect to the lexicographic order.
Then we define

\begin{equation}
\label{eq:action_aux}
\tilde A^n_\omega(\gamma)=\sum_{k=0}^{n-1}L(\tilde\gamma_{k+1}-\tilde \gamma_k)+\sum_{k=0}^{n-1} F_\omega(k,\tilde\gamma_k).
\end{equation}
and 
\[
\tilde A^n_\omega(x)=\inf\big\{\tilde A^n_\omega(\gamma):\gamma\in\Gamma^n_\omega,\ \gamma_n=x\big\}.
\]
Although the concatenation of two admissible paths is not admissible, replacing the non-Poissonian point with its Poissonian approximation we obtain an admissible path with action equal to the sum of actions of individual paths. 
Thus the subadditivity property~\eqref{eq:subadditivity} holds 
almost surely for $\tilde A^n_\omega$.

To apply the subadditive ergodic theorem, we need upper and lower estimates on $\tilde A^n_\omega(nv)$. The lower linear estimate
$\tilde A^n_\omega(nv)\ge bn$
follows from~\eqref{eq: growth of L}.

To prove an upper estimate, let's fix any $v\in\R^d$. For all $k=0,1,2,\ldots$, we define $\tilde\gamma_k$ as the spatial component $x$ of the Poissonian point $(k,x)$ closest to $(k,kv)$ (using the lexicographic order to break ties). Then
\[
\tilde A^n_\omega(nv)\le \sum_{k=0}^{n-1}L(\tilde\gamma_{k+1}-\tilde\gamma_{k})+ \sum_{k=0}^{n-1} F_\omega(k,\tilde\gamma_k),
\]
and to obtain a linear upper bound,
by stationarity and the ergodic theorem, it suffices to check that the expectation of each term on the right-hand side is finite.
For the terms in the second sum, it follows from our assumption~\eqref{eq:labels are L1}. 
The terms in the first sum form a stationary ergodic sequence. To apply the ergodic theorem, we will derive the finiteness of expectation from~\eqref{eq:comparing to ball volume} and~\eqref{eq:assumption on growth of L}.

 First, recall that for every $\eta>0$ and all $x,y>0$,
\begin{equation}
\label{eq:elementary}
(x+y)^d\le (1+\eta)^dx^d+K_\eta y^d,
\end{equation}
where $K_\eta=(1+\eta^{-1})^d$. To see this, consider the cases where $x<y\eta^{-1}$ and $x\ge y \eta^{-1}$.

Let's denote $\xi_k=\tilde\gamma_k-kv$.
Using~\eqref{eq:assumption on growth of L} and the triangle inequality,  plugging  $x=|\xi_{k+1}|+|\xi_{k}|$ and $y=|v|$ and into~\eqref{eq:elementary},  we
obtain
\begin{equation}
\label{eq:moment computation}
\begin{split}
    \E L(\tilde\gamma_{k+1}-\tilde\gamma_{k})
    &=\E L(v+\xi_{k+1}-\xi_k) 
    \\
    &\le C_1 e^{C_2 K_\eta|v|^d} \E e^{C_2(1+\eta)^d(|\xi_{k+1}|+|\xi_k|)^d}
    \\
    & \le C_1 e^{C_2 K_\eta|v|^d} \E e^{C_22^{d-1}(1+\eta)^d|\xi_{k+1}|^d} \E e^{C_22^{d-1}(1+\eta)^d|\xi_k|^d},
\end{split} 
\end{equation}
where in the second line we used $(x+y)^d\le 2^{d-1}(x^d+y^d)$ valid for all $x,y>0$.
Since for all $k$,
\begin{equation}
\label{eq:tail of distance to closest}
\Pp\{|\xi_k|>z\}=e^{-V_dz^d},\quad z>0, 
\end{equation}
we can use \eqref{eq:comparing to ball volume} to
obtain that the right-hand side of~\eqref{eq:moment computation} is finite if $\eta>0$ is chosen sufficiently small. 
This concludes the proof of at most linear growth for~$\tilde A^n_\omega(nv)$.
Now the subadditive ergodic theorem implies that there is a function $\Lambda$ such that for each $v\in\R^d$,
\begin{equation}
\label{eq:shape for aux}
    \lim_{n\to\infty}\frac{\tilde A^n_\omega(nv)}{n}\aseq\Lambda (v).
\end{equation}

The standard argument based on the fact that the concatenation of two paths of different slopes provides  an upper estimate for the minimal action for the averaged slope yields convexity of $\Lambda$.

\smallskip

Now let's prove that $\Lambda$ obtained in~\eqref{eq:shape for aux} also serves as the shape function for the original action $A^n$. 
We need to prove an upper and a lower bound on $A^n_\omega(nv)$. For an upper bound, we take any $v$, and  write
\begin{equation}
\label{eq:upper estimate real vs aux}
A^n_\omega(nv)\le L(\tilde\gamma_1-0)+\tilde A^{n-2}_{\Theta^{1,v}\omega}((n-2)v)+L(nv-\tilde\gamma_{n-1}),\quad n\ge 2,
\end{equation}
where $\tilde\gamma_k$ has been defined above. 
The first term on the right-hand side does not depend on $n$, the second term satisfies \eqref{eq:shape for aux} with $\omega$ replaced by
its shift $\Theta^{1,v}\omega$ since $\Theta^{1,v}$ preserves $\Pp$. For the third term, we note that
the sequence of r.v.'s $(L(nv-\tilde\gamma_{n-1}))_{n\in\N}$ is stationary and ergodic, and since a computation similar to \eqref{eq:moment computation}
implies that they have finite expectation, we conclude that $L(nv-\tilde\gamma_{n-1})=o(n)$ holds with probability $1$. Combining this with~\eqref{eq:shape for aux} and~\eqref{eq:upper estimate real vs aux}, we obtain
\begin{equation}
\label{eq:upper bound wrt aux action}
\limsup_{n\to\infty}\frac{A^n_\omega(nv)}{n}\le \Lambda (v)\ \mathrm{a.s.}
\end{equation}

Now let's obtain a lower bound.
Let $\gamma=\gamma^n_{A,\omega}(nv)$ be the minimizer realizing the minimum in the definition of $A^n_\omega(nv)$ (see Proposition~\ref{prop:min}). Replacing the endpoints of $\gamma$ by Poissonian points $(0,\tilde\gamma_0)$ and $(n,\tilde\gamma_n)$ we can bound the modified action with
\begin{multline}
\label{eq:estimating from below}
\tilde A^n_\omega(nv)\le F_\omega(0,\tilde\gamma_0)+A^n_\omega(nv)\\ +[L(\gamma_1-\tilde\gamma_0)-L(\gamma_1-0)]
+[L(\tilde \gamma_n-\gamma_{n-1})-L(nv-\gamma_{n-1})].
    \end{multline}

Since 
\[A^n_\omega(nv)=A^n_\omega(\gamma)\ge \sum_{k=0}^{n-1}L(\gamma_{k+1}-\gamma_k),   
\]
and $L(\gamma_{k+1}-\gamma_k)\ge b$ for all $k$, we obtain
\[
\max\{L(\gamma_1), L(nv-\gamma_{n-1})\}\le A^n_\omega(nv)- (n-1)b.
\]
The event  where $A^n_\omega(nv)\le (\Lambda(v)+1)n$ for all large $n$ has probability~$1$ by~\eqref{eq:upper bound wrt aux action}.
On this event, for those large $n$, we have
\[
\max\{L(\gamma_1),L(nv-\gamma_{n-1})\}\le s_n,
\]
where $s_n=(\Lambda(v)+1-b)n+|b|$. It is easy to see that $\Lambda(v)\ge b$, so $s_n\ge n$. 

To apply~\eqref{eq:cond on large increments of L} to \eqref{eq:estimating from below}, we also need an estimate on
the increment of the argument for the last two terms. For the first of them, 
\[
|(\gamma_1-\tilde\gamma_0)-(\gamma_1-0)|=|\tilde \gamma_0|
\]
does not depend on $n$. The second one satisfies
\[
|(\tilde \gamma_n-\gamma_{n-1})-(nv-\gamma_{n-1})|=|\tilde \gamma_n-nv|=|\xi_n|< c (\ln n)^{1/d}\le c (\ln s_n)^{1/d}
\]
for sufficiently large $n$ due to \eqref{eq:tail of distance to closest} and the Borel---Cantelli lemma, if \eqref{eq:parameter c} holds.
Thus, both of these argument increments are bounded by $c(\ln s_n)^{1/d}$, and due to~\eqref{eq:cond on large increments of L}, \eqref{eq:rho},
and~\eqref{eq:estimating from below},
we obtain that with probability~$1$,
\[
A^n_\omega(nv)\ge \tilde A^n_\omega(nv) - F_\omega(0,\tilde\gamma_0)-2\rho_c(s_n)=\tilde A^n_\omega(nv) - F_\omega(0,\tilde\gamma_0)+o(n), 
\]
which, together with 
\eqref{eq:shape for aux} implies
\[
\liminf_{n\to\infty}\frac{A^n_\omega(nv)}{n}\ge \Lambda(v),
\]
completing the proof.

\providecommand{\bysame}{\leavevmode\hbox to3em{\hrulefill}\thinspace}
\providecommand{\MR}{\relax\ifhmode\unskip\space\fi MR }
\providecommand{\MRhref}[2]{%
  \href{http://www.ams.org/mathscinet-getitem?mr=#1}{#2}
}
\providecommand{\href}[2]{#2}

\end{document}